\documentclass[10pt]{article}
\usepackage[latin1]{inputenc}
\usepackage[english]{babel}
\usepackage{amssymb,amsfonts, amsmath}
\usepackage{graphicx, xcolor}
\newtheorem{theorem}{Theorem}[section]
\newtheorem{proposition}[theorem]{Proposition}
\newtheorem{corollary}[theorem]{Corollary}
\newtheorem{lemma}[theorem]{Lemma}
\newtheorem{remark}[theorem]{Remark}

\newcommand{\bcl}{\begin{center}}
\newcommand{\ecl}{\end{center}}
\newcommand{\brl}{\begin{right}}
\newcommand{\erl}{\end{right}}
\newcommand{\ben}{\begin{enumerate}}
\newcommand{\een}{\end{enumerate}}
\newcommand{\barr}{\begin{array}}
\newcommand{\earr}{\end{array}}
\newcommand{\btab}{\begin{tabular}}
\newcommand{\etab}{\end{tabular}}
\newcommand{\bdoc}{\begin{document}}
\newcommand{\edoc}{\end{document}}
\newcommand{\beqy}{\begin{eqnarray}}
\newcommand{\eeqy}{\end{eqnarray}}

\newcommand{\beqi}{\begin{eqnarray*}}
\newcommand{\eeqi}{\end{eqnarray*}}
\newcommand{\bitem}{\begin{itemize}}

\newcommand{\eitem}{\end{itemize}}
\newcommand{\nln}{\newline}
\newcommand{\newt}{\newtheorem}
\newcommand{\pa}{\partial}
\newcommand{\re}{{I\!\!R}}
\newcommand{\ren}{\re^N}
\newcommand{\xr}{x\in\re }
\newcommand{\x}{\times}
\newcommand{\dyle}{\displaystyle}
\newcommand{\ene}{{I\!\!N}}
\newcommand{\irn}{\int\limits_{\re^N}}
\newcommand{\io}{\int\limits_{\O}}
\newcommand{\meas}{{\rm meas\,}}

\newcommand{\sign}{{\rm sign}}
\newcommand{\map}{\longrightarrow }
\newcommand{\imp}{\Longrightarrow }
\newcommand{\sen}{{\rm sen\,}}
\newcommand{\tg}{{\rm tg\,}}
\newcommand{\arcsen}{{\rm arcsen\,}}
\newcommand{\arctg}{{\rm arctg\,}}
\newcommand{\supp}{{\textsl supp\ }}
\newcommand{\ity}{\int_{-\iy}^{+\iy}}
\newcommand{\limit}{\lim\limits}
\newcommand{\limi}{\limit_{n\to\infty}}
\newcommand{\sumi}{\sum\limits_{n=1}^{\infty}}
\newcommand{\ulu}{\underline u}
\newcommand{\ulw}{\underline w}
\newcommand{\ulz}{\underline z}
\newcommand{\ulv}{\underline v}
\newcommand{\uls}{\underline s}
\newcommand{\olu}{\overline u}
\newcommand{\olv}{\overline v}
\newcommand{\ols}{\overline s}
\newcommand{\ob}{\overline\b}
\newcommand{\ovar}{\overline\var}
\newcommand{\wv}{\widetilde v}
\newcommand{\wu}{\widetilde u}
\newcommand{\ws}{\widetilde s}
\renewcommand{\a }{\alpha }
\renewcommand{\b }{\beta }
\newcommand{\g }{\gamma}
\newcommand{\G }{\Gamma }
\renewcommand{\d }{\delta }

\newcommand{\D }{\Delta }
\newcommand{\e }{\varepsilon }
\newcommand{\z }{\zeta }
\renewcommand{\l }{\lambda }
\renewcommand{\L }{\Lambda }
\newcommand{\m }{\mu }
\newcommand{\n }{\nabla }
\newcommand{\s }{\sigma }
\newcommand{\Sig }{\Sigma }
\renewcommand{\t }{\tau }
\newcommand{\var }{\varphi }
\renewcommand{\o }{\omega }
\renewcommand{\O }{\Omega }
\newcommand{\bR}{{\bf R}}
\newcommand{\bC}{{\bf C}}
\newcommand{\bZ}{{\bf Z}}
\newcommand{\bN}{{\bf N}}
\newcommand{\bQ}{{\bf Q}}
\newcommand{\bK}{{\bf K}}
\newcommand{\bI}{{\bf I}}
\newcommand{\bv}{{\bf v}}
\newcommand{\bV}{{\bf V}}
\DeclareMathOperator{\di}{div} \DeclareMathOperator{\Hess}{Hess}
\DeclareMathOperator{\conv}{conv} \DeclareMathOperator{\cd}{D}
\DeclareMathOperator{\ud}{\overline{D}}
\DeclareMathOperator{\suppo}{supp} \DeclareMathOperator{\Ric}{Ric}
\DeclareMathOperator{\del}{\delta} \DeclareMathOperator{\inj}{inj}
\def\qed{\unskip\kern 6pt \penalty 500
\raise -2pt\hbox{\vrule \vbox to10pt{\hrule width 4pt
\vfill\hrule}\vrule}\par}

\newenvironment{Proof}{\removelastskip\vskip12pt
plus 1pt \noindent\em\rm}{\hfill {\qed \hskip .2cm}}

\title
{Allen-Cahn Approximation of Mean Curvature Flow in Riemannian
manifolds II, \\Brakke's flows}

\author{Adriano Pisante\thanks{Dipartimento di Matematica "G.
Castelnuovo", Universit\`a di Roma ``La Sapienza'', P.le A.~Moro
5, I-00185 Roma, Italia (pisante@mat.uniroma1.it).} \ \ and Fabio
Punzo\thanks{Dipartimento di Matematica "G. Castelnuovo",
Universit\`a di Roma ``La Sapienza'', P.le A.~Moro 5, I-00185
Roma, Italia (punzo@mat.uniroma1.it).}}

\date{}

\begin{document}
 \maketitle
{\abstract{ \noindent We prove convergence of solutions to the
parabolic Allen-Cahn equation to Brakke's motion by mean curvature
in Riemannian manifolds, generalizing previous results from \cite{Ilm1} in
Euclidean space. We show that a sequence of measures, associated
to energy density of solutions of the parabolic Allen-Cahn
equation, converges in the limit to a family of rectifiable Radon measures, which
evolves by mean curvature flow in the sense of Brakke. A key role
is played by a local almost monotonicity formula (a weak
counterpart of Huisken's monotonicity formula) proved in
\cite{PiPu1}, to get various density bounds for the limiting
measures.

\bigskip

%\noindent {\it  2010 Mathematics Subject Classification: .}

}}

\bigskip
\medskip
\smallskip

\section{Introduction} \setcounter{equation}{0}
In \cite{PiPu1} we started to investigate the Allen-Cahn equation
\begin{equation}\label{e1}
\pa_t u^\e \,=\, \Delta u^\e -\frac 1{\e^2} f(u^\e)\quad
\textrm{in}\;\; M\times(0,\infty)\,,
\end{equation}
completed with the initial condition
\begin{equation}\label{e1a}
u^\e=u^\e_0\quad \textrm{in}\;\; M\times \{0\}\,,
\end{equation}
were $\e>0$ is a small parameter and $M$ is an $N-$dimensional
Riemannian manifold with Ricci curvature bounded from below.

We suppose that the nonlinearity is the negative gradient of a double well potential $F$ with two minima of equal depth. More precisely, we assume for simplicity that
\[
\textrm{\ \ } \left\{
\begin{array}{l}
(i) \quad\,\; f = F',\; \textrm{with}\;\; F\in C^\infty(\re),\,
F\;\,\textrm{even}\,;
\\
(ii) \quad f(0)=f(\pm1)=0\,, f<0\,\,\hbox{in}\,\, (0,1)\,, f>0\,\,\hbox{in}\,\, (1,\infty),\,\\ \qquad \;\, f'(0)<0,
f'(\pm1)>0\,;
\\ (iii)\;\; \,  F>0\;\,\textrm{in}\;\, \re\setminus\{\pm1\},\; F(\pm1)=0\,;
\\ (iv) \;\;\, \min_{[\a, \infty)}F''>0,\, \textrm{for some}\;\, \a\in
(0,1)\,.
\end{array}
\right.\leqno(H_0)\]

For any $\e>0$ and for any $(x,t)\in M\times [0,\infty)$, set
\begin{equation}\label{e2}
E^\e(x,t):=\frac 1 2 |\nabla u^\e|^2+\frac 1 {\e^2} F(u^\e)\quad
(x\in M, t\ge0)\,,
\end{equation}
and define the energy density
\begin{equation}\label{e8i}
d\mu^\e_t := \left\{\frac{\e}2 |\nabla u^\e|^2 + \frac
1{\e}F(u^\e)\right\}d\mathcal V(x).
\end{equation}
Clearly,
\begin{equation}\label{e3}
d\mu^\e_t(x)=\e E^\e(x,t) d \mathcal V(x)\quad (x\in M, t\ge 0)\,.
\end{equation}
\smallskip

Let $E_0\subset M$ be an open bounded subset with smooth boundary
$\pa E_0=\Sigma_0$. Hence there exist $C_0>0, R_0>0$ such that
$$\mathcal H^{N-1}\big(\Sigma_0\cap B_R(x)\big)\le C_0
\omega_{N-1}R^{N-1}$$ for all $0<R<R_0$\,. In the end, the hypersurface $\Sigma_0$ will be the initial datum for an evolution by mean curvature, in a suitable weak sense, obtained as limit of diffuse interfaces, the regions say $\{ |u^\e|\leq \alpha \}$, where $u^\e$ solve \eqref{e1}-\eqref{e1a}.

Concerning the initial conditions $u_0^\e$ in \eqref{e1a} (and the corresponding
$\mu^\e_0\equiv \mu^\e(\cdot,0)$ given by \eqref{e8i}) we always
assume the following:
\[\left\{
\begin{array}{l}
(i)\,\,\mu^\e_0\to \a  \mathcal H^{N-1}\lfloor
\Sigma_0\,\;\hbox{as}\,\, \e\to 0\,\,\hbox{as Radon measures, for
some}\,\,\a\ge 0\,;
\vspace{.2 cm} \\
(ii)\,u_0^\e \to 2 \chi_{E_0} -1\quad \textrm{as}\;\,\e\to\infty\,\,*-\textrm{weakly in}\;\, BV_{loc}(M)\, ; \\
(iii)\,\hbox{there exists}\,\, C_0>0\,\,\hbox{such that}\,\,
\frac{\mu_0^\e(B_R(x))}{\omega_{N-1}R^{N-1}}\le C_0 \\ \quad
\;\;\;\;
\textrm{for all}\; x\in M,\, 0<R<R_0,\, 0<\e<1; \vspace{.2 cm} \\
(iv)\, \hbox{there exists}\,\, k_0>0\,\, \hbox{such that}\,\,
\|u_0^\epsilon\|_\infty\leq H_0\,;\\
(v)\,\; \hbox{there exists}\,\, \check C>0\,\,\hbox{such that for
any}\,\, 0<\e<1 \\ \quad \;\,\; \|\nabla u_0^\e\|_\infty\le
\frac{\check C}{\e} \,. \vspace{.2 cm} \\
\end{array}
\right. \leqno(H_1)
\]
Under hypotheses $(H_0), (H_1)$ in \cite{PiPu1} it is shown that problem
\eqref{e1}-\eqref{e1a} admits a unique bounded solution. Moreover,
$u^\e\in C^{\infty}\big(M\times(0,\infty)\big)\cap
C\big(M\times[0,\infty)\big);$
\begin{equation}\label{e121b}
|u^\e | \leq  k_0 \quad \textrm{for all}\;\; x\in M, t>0\,.
\end{equation}
In addition,
\begin{equation}\label{e162}
\sup_{\e>0} \sup_{t\in (0,\infty)} \e \int_M E^\e d\mathcal V \le
C_2\,
\end{equation}
where $C_2:=\sup_{\e>0}\mu_0^\e(M)$;
\begin{equation}\label{e150}
t\mapsto \int_M E^\e(x,t) d\mathcal V(x)\,\;\,\textrm{is
nonincreasing for}\;\; t>0\,.
\end{equation}

Moreover, recalling the difinition of discrepancy Radon measure
\[d \xi_t^\e:=\left(\frac{\e}2 |\nabla u^\e|^2-\frac 1{\e} F(u^\e)\right)d \mathcal V\,,\]
it is proved that
\begin{equation}\label{u13}
\limsup_{\e\to 0^+} \sup_{(x,t)\in Q} \xi^\e_t(x)\leq 0\,,
\end{equation}
for each compact subset $Q\subset M\times (0,\infty)$. Using
\eqref{u13} is then shown that
\begin{equation}\label{e60}
\begin{split}
\frac{d}{dt}\int_{M} \phi(x,t)d\mu^\e_t\,\le\,
\frac{C_3}{\sqrt{s-t}}\int_{M} \phi(x,t)d\mu^\e_t + C_4 + \frac
{C_5}{\sqrt{s-t}}
\end{split}
\end{equation}
for all $0\le t<s$, for some positive constants $C_3, C_4, C_5$
independent of $\e$.

Here, for any fixed reference point $(y,s)\in
M\times (0,\infty)$, $\phi(x,t)\equiv \phi(x,t; y,s)$ is a
suitable kernel, depending explicitely on the
Riemannian distance $d(x)=d(x,y)$ for $x,y\in M$ as follows
\begin{equation}
\label{kernel}
\phi(x,t)=\hat{\zeta}(d^2(x)) (s-t)^{-\frac{N-1}{2}} e^{-\frac{d^2(x)}{4(s-t)}} \, .
\end{equation}
 In constrast with the case of $\re^N$, it has a suitably small
compact support in space due to the cut-off function $\hat{\zeta}$. In addition, we were able to deduce the previous inequality in full generality, without assuming well prepared initial data as in \cite{Ilm1}.

 Clearly, inequality \eqref{e60} does not
imply monotonicity for the function $t\to\int_{M}
\phi(x,t)d\mu^\e_t(x)$. Nevertheless,
it still allows us to
control the behavior of $d\mu^\e_t$ at small scales. For this
reasons, we refer to \eqref{e60} as a {\it local almost
monotonicity formula}, in analogy with the monotonicity formula valid in the Euclidean space.

As a consequence of \eqref{e60}, we were able to obtain the inequality
\begin{equation}\label{e61}
\begin{split}
\mathcal G(t)\le e^{\frac{C_3}2 (\sqrt{s-t_0}  -\sqrt{s-t})}\big[ \mathcal
G(t_0) + C_4(t-t_0) + C_5( \sqrt{s-t_0}  -\sqrt{s-t}) \big]\,,
\end{split}
\end{equation}
 for all $0\leq
t_0<t<s$, where
\[\mathcal G(t):= \int_{M} \phi(x,t) d\mu^\e_t\quad (0\le t<s)\,.\]

Actually, \eqref{u13}, \eqref{e60} and \eqref{e61} remain true (see
\cite{PiPu1} for details), if instead of $(H_1)$, we only assume that
\eqref{e121b} is satisfied, and that for each compact subset
$K\subset M, T>0$, there holds:
\begin{equation}\label{e170}
\sup_{\e>0} \sup_{t\in (0,T)} \e \int_K E^\e d\mathcal V \le
\overline C\,
\end{equation}
for some constant $\overline C>0$ depending on the compact subset
$K$ and on $T$, but independent of $\e$.
In addition, for each compact subset $K\subset M$ we have that
\[\int_{M}\phi(x,t; y,s) d\mu^\e_t(x)\le \underline C\quad \textrm{for all}\;\; y\in M,  0\le t<s\,, \leqno (G^\e_1)\]
for some $\underline C=\underline C_K>0$, and
\[\mu^\e_t(B_R(x))\le \omega_{N-1} D_0 R^{N-1}\quad \textrm{for all}\;\; x\in K, 0< R<\tilde R, t\ge 0, \leqno (G^\e_2)\]
for some $0<\tilde R<R_0$ and $D_0=D_0(\underline C, \tilde R)>0$
independent of $\e$\,.

\smallskip

In the present paper, extending results from \cite{Ilm1} in the Euclidean space, we describe the asymptotic behaviour of  the
family of measures $\{\mu^\e_t\}$ as $\e \to 0$ on Riemannian manifolds.
First of all, adapting the semidecreasing trick from \cite{Ilm1}-\cite{Ilm2}, we prove that we can extract a subsequence
$\mu^{\e_n}_t$ that, for every $t>0$, converges as Radon measure
on $M$ to a limit Radon measure $\mu_t$ for all $t\ge 0$ as $n\to
\infty$.

%Then, exploiting the structure of space form we show the {\it Empty Spot Lemma} (see Subsection \ref{esl}), a technical result which, under suitable assumptions, ensures locally the absence of interface in a pointwise sense. Actually this is the only point of the paper where we need that
%$M$ is a space form. Although we firmly believe that such property
%holds in full generality, extending the corresponding results we use from the Euclidean space (see \cite{Chen}; see also \cite {AHM}) to the case of arbitrary manifolds seems at present highly non obvious.

Using \eqref{e60}, a
version of Brakke's {\it Clearing-Out Lemma} \cite[Lemma
6.3]{Brakke} is established, similar to \cite[Lemma 6.1]{Ilm1}. Note that our proof is direct and self-contained. In particular, it does not rely on the so-called {\it Empty Spot Lemma} \cite[Lemma 6.4]{Ilm1} and on the related results on propagation of fronts (see \cite{Chen}), which at present are not available on Riemannian manifolds. However, we also point out that  the same strategy as \cite{Ilm1} could be adapted to the present situation (see Remark \ref{emptysp}). Furthermore, we do not use the gradient bound coming from the assumption of well prepared data. In particular, it is shown that if
$(y,s)\not\in \overline{\bigcup_{t'\ge 0} \suppo
\mu_{t'}\times\{t'\}}$, then there exists a neighborhood $U\subset
M\times[0,\infty)$ of $(y,s)$ such that $\{u^{\e_n}\}$ converges
uniformly in $U$ to either $1$ or $-1$ as $n\to \infty$\,. Such a result, well known in the Euclidean case, show once more absence of evolving interface where there is no energy concentrating in the limiting measures $\mu_t$. In addition, an estimate for the size of the bad set follows, showing that $\mathcal H^{N-1}(\suppo \mu_t)$ is locally finite for a.e. $t>0$.

With the local almost monotonicity formula \eqref{e60} at disposal, we
adapt to the present situation the strategy of \cite{Ilm1} and show that
the discrepancy measure $d\xi^\e_t$ converges to
$0$ as $\e \to 0^+$. Indeed, in view of \eqref{u13} it is enough to consider the negative part of the limiting discrepancy $d\xi_t$. In addition, at $|d\xi|-$a.e. point in space-time a suitable (forward) density of $\mu_t$ defined through \eqref{kernel} is  shown to be both zero (as a consequence of \eqref{e60}) and strictly positive (because of the Clearing-Out lemma), so the discrepancy has to vanish identically.

Thus, we obtain all preliminary results necessary to pass to limit as
$\e\to 0^+$, in the sense of varifolds, in the Brakke's type
equality \eqref{e175} satisfied by $\mu^\e_t$ (see Section
\ref{vld}). Hence $N-1-$rectifiability for the limiting measures $\mu_t$ for a.e. $t>0$ and Brakke's inequality, namely

\begin{equation}\label{e90}
\ud_t \mu_t(\phi)\, \le\, \int_{M}\big\{ -\phi H_t^2 + \langle \nabla \phi,
T_x\mu_t^{\perp}(\,\overrightarrow{H}_t\,)\rangle\,\big\}\,d\mu_t\,,\,
\end{equation}
for all $ \phi\in C^2_c(M ; \re^+)\,$ and for every $t>0$, follows at once (all the terms in the formula being actually well defined and finite for a.e. $t>0$). Here $\ud_{t}$ is the {\it upper derivative} of
\[\mu_t(\phi)\equiv \int_{M} \phi(x) d\mu_t(x)\, , \]
$\overrightarrow{H}_t$ is the mean curvature vector associated to the varifold corresponding to $\mu_t$  and $T_x \mu_t^\perp$ is the orthogonal projection onto the normal space to the measure  (see Section 2 for precise definitions). Note that in this paper we do not address the issue of integrality for the limiting rectifiable measures $\mu_t$ (or, equivalently, for the corresponding varifolds). This should follow from a careful adaptation of the subtle result in \cite{Ton1} (see also \cite{Ton10}) valid in the Euclidean space. In particular, once integrality is established one would have $\overrightarrow{H} \perp T_x \mu_t$ a.e. (see \cite{Brakke}) and in particular no projection operator in \eqref{e90}. In addition, we do not investigate partial regularity property  of the solution we contruct. In this respect, when trying to discuss this issue, expecially in connection with the so-called 'unit density hypotesis' for the limiting varifolds, it would be natural to generalize the recent partial regularity results from \cite{KT}, \cite{Ton2} valid in the Euclidean case.

 As a final remark we observe that, among several possible ways to obtain global weak solutions of the mean curvature flow on manifolds, such as the level-set approach via viscosity solutions (see e.g. \cite{ES2}, \cite{Ilm3}, \cite{AAM} or \cite{CGG}), the method of barriers (see e.g \cite{BellNov},
\cite{BellPaol})) or the geometric measure
theory approach via varifolds, currents or BV functions (see, e.g. \cite{Brakke}, \cite{Ilm2} or \cite{LS}), we decided to adapt the Allen-Cahn approximation from \cite{Ilm3}. In our opinion this approach seems more promising in order to flow unbounded initial hypersurfaces with only locally finite area. Such problem arises naturally for example when trying to evolve complete noncompact surfaces in the hyperbolic space. Indeed, unbounded minimal hypersurfaces with prescribed boundary at infinity in $\mathbb{H}^N$ exist in abundance and can be constructed e.g. by the stationary phase-field approximation analogous to \eqref{e1} (see e.g. \cite{PisPons}). As it was the main motivation for the present research, we plan in a future paper to study convergence to such equilibria under mean curvature flow in $\mathbb{H}^N$ for unbounded hypersurfaces with fixed boundary at infinity and the connections of such evolution with the renormalized area studied e.g. in \cite{AM}.

\section{Mathematical background: varifolds, rectifiable Radon measures, rectifiable varifolds, first
variation, mean curvature}\label{varif} In this Section we
recall some preliminaries from Geometric Measure Theory (for more
details see, $e.g.$, \cite{All}, \cite{Sim}, \cite{Ilm2}).

To begin with, recall that by Nash Embedding Theorem, we can
assume that $M$ is isometrically embedded in $\re^L$ for some
$L\ge N$.

Let $G(L,k)$ be the Grassman manifold of unoriented $k-$planes in
$\re^L\;\; (k\leq N)$; let $$G_k(M):=\big\{(x,S)\in M\times
G(L,k)\,:\, S \subset T_x M \big\}\,.$$

A {\it general} $k-$ {\it varifold} is a Radon measure on
$G_k(M)$. We denote by $\mathbb V_k(M)$ the set of all general
$k-$varifolds, and we give it the topology corresponding to the weak-$*$ convergence
of Radon measures. We write:
\[V(\psi)=\int_{G_k(M)}\psi (x,S) dV(x,S),\]
where $\psi\in C^0_c(G_k(M)).$

\smallskip

Denote by  $\mathcal M^+(M)$ the set of all positive Radon
measures on $M$. Given $\mu\in \mathcal M^+(M)$ For any $x\in M,
\l>0, 1\le k\le N$ we can define the scaled Radon measure
$\mu_{x,\l}\in \mathcal M(\re^L)$ by
\[\mu_{x,\l}(A)=\frac{\mu [M\cap (\l A+x )]}{\l^k}\quad
(A\subseteq \re^L)\,.\] Let $\mathcal P$ a $k-$plane in $T_xM$ and
$\a>0$. We say that $\mathcal P \equiv T_x\mu $ is the
$k-$dimensional {\it approximate tangent plane} of $\mu$ at $x$,
if
\[\lim_{\l\to 0^+}\mu_{x,\l}=\a \mathcal H^k\lfloor \mathcal P\,,\]
where $\mathcal H^k$ is the $k-$dimensional Hausdorff measure in
$\re^L$. The $k-$dimensional approximate tangent space of a set
$X\subseteq M$ at $x\in M$ is defined by
\[T_x X:=T_x(\mathcal H^k\lfloor X)\,,\]
if it exists. We say that $X$ is {\it countably} $k-${\it
rectifiable}, if $X\subseteq C_0\cup (\cup_{i\ge 1} C_i)$, where
$\mathcal {H}^k(C_0)=0$ and each $C_i\subseteq M$ is an embedded
$C^1\, k-$ submanifold. We say that $X\subseteq M$ is {\it
locally} $k-${\it rectifiable} if in addition $X$ has locally
finite $\mathcal {H}^k-$measure. If $X$ is locally $k-$rectifiable
and $\mathcal{H}^k-$measurable, then $T_x X$ exists $\mathcal
H^k\lfloor X-$a.e. .

\smallskip

Let $X$ be an $\mathcal H^k-$measurable subset of $M$, and let
$\theta:M\to [0,\infty)$ be locally $\mathcal H^k-$integrable,
with $X=\{\theta>0\}\,\mathcal H^k-$a.e. . The Radon measure $\mu$
on $M$ is a $k-${\it rectifiable Radon measure}, if either of the
following equivalent conditions holds:
\begin{itemize}
\item[$(a)$] $\mu$ has $k-$dimensional tangent planes of positive multiplicity $\mu-$a.e.;
\item[$(b)$] $\mu=\mathcal H^k\lfloor \theta$ for some $X$ which is $\,\mathcal
H^k-$measurable and countably $k-$rectifiable, and $\theta$
locally $\mathcal H^k-$integrable.
\end{itemize}
We denote by $\mathcal M_k(M)$ the set of $k-$rectifiable Radon
measures on $M$. We call $\mu$ an {\it integer} $k-${\it
rectifiable Radon measure}, if either of the following equivalent
conditions hold:
\begin{itemize}
\item[$(c)$] $\mu$ has $k-$dimensional tangent planes of positive integer multiplicity $\mu-$a.e.;
\item[$(d)$] $\mu=\mathcal H^k\lfloor \theta$ for $X\,\mathcal
H^k-$measurable and locally $\mathcal H^k-$rectifiable, and
$\theta$ locally $\mathcal H^k-$integrable with values in $\ene$.
\end{itemize}
We write $\mathcal I\mathcal M_k(M)$ for the set of all such
$\mu$.

\smallskip

Associated to a varifold $V$ there is a Radon measure $\mu_V$ on
$M$, defined by
\[\mu_V:=\pi_{\sharp}(V)\,,\]
where $\pi:G_k(M)\to M$ is the natural projection map. On the other hand,
if $\mu$ is a $k-$rectifiable Radon measure, then $x\mapsto
T_x\mu$ is a $\mu-$measurable section of $G_k(M)$ defined
$\mu-$a.e.\, . Therefore, we can define the varifold $V=V_\mu$ by
\[V_\mu(\psi)=\int_{G_k(M)}\psi(x,S) dV_\mu(x,S)=\int_{M}\psi(x,T_x\mu)d\mu(x)\]
for $\psi\in C^0_c(G_k(M))\,.$ Observe that $\mu= \|V\|,$ where
$\|V\|$ is the Radon measure defined as follows:
\[\|V\|(A):=V\Big(G_k(M)\cap\pi^{-1}(A\Big)\quad \textrm{whenever}\;\; A\subset M\,.\]

We denote by $\mathcal R \mathbb V_k(M)$ the set of $k-${\it
rectifiable varifolds}, $i.e.$ the varifolds associated to
$k-$rectifiable Radon measures on $M$; whereas, by $\mathcal
I\mathbb V_k(M)$ the set of {\it integer} $k-${\it rectifiable
varifolds}, $i.e.$ the varifolds associated to integer
$k-$rectifiable Radon measures on $M$.

\smallskip

When $S$ is a $k-$plane in $G(T_x M, k)$, we also use $S$ to
denote the orthogonal projection from $T_xM$ onto $S$.
Furthermore, we write $A:B$ for the inner product of endomorphisms
$A$ and $B$ on $T_xM$.

Let $U\subseteq M$ be an open subset; we set $C^1_c(U,
TU):=\big\{Y\in \Gamma(TM)\,:\,\suppo Y\subset U \big\}\,.$ When
$U=M$, we write $C^1_c(TM)$ instead of $C^1_c(M, TM)$\,.

Let $Y\in C^1_c(TM).$ Recall the {\it first variation formula} for
a varifold
\[\del V(Y)=\int_{G_k(M)} S : DY dV(x,S)\,;\]
here
$$\big(S:\cd Y\big)(x)=\sum_{i=1}^k\langle \cd_{e_i} Y(x), e_i \rangle, $$
$\{e_1, \ldots, e_k\}$ being any orthonormal basis of $S$ and $\cd_{e_i}  Y$ the corresponding covariant derivatives. For $U\subseteq M$ open, define the {\it total first variation} by
\[|\del V|(U)=\sup \big\{\del V(Y)\,:\,Y\in C^1_c(U), |Y|\le 1 \big\}\,.\]
If $\d V$ is a Radon measure and $|\d V|<<\|V\|$, then
\[\del V(Y)=\int_{G^k(M)} S : DY dV(x,S)=-\int_{M}\langle Y(x), \overrightarrow{H}(x)\rangle d\|V\|(x),\]
where $\overrightarrow{H}: M \to TM$ is the Radon-Nikodym
derivative of $\delta V$ with respect to $\| V\|$. By definition
$\overrightarrow{H}$ is called the {\it mean curvature vector
field}. When, in addition, $V=V_\mu$ is also a $k-$rectifiable
varifold, then
\[\del V_\mu(Y)=\int_{M}T_x\mu :DY d\mu=-\int_{M}\langle\overrightarrow{H},Y\rangle d\mu\,.\]

\smallskip

%\medskip

%We shall use the following celebrated Allard compactness theorem
%for varifolds (see \cite{All}, \cite{Ilm2}).
%\begin{theorem}\label{thmAll}
%Let $\{\mu_n\}\subset \mathcal I \mathcal M_k(M)$ with
%\[\sup_{n\in \ene}\big[\mu_n(U)+|\del V_{\mu_n}|(U)\big] < \infty\quad \textrm{for each}\;\; U \subset\subset M\,.\]
%Then there exists $\mu\in \mathcal I\mathcal M_k(M)$ and a
%subsequence $\{\mu_{n_k}\}\subset \{\mu_n\}$ such that
%\begin{itemize}
%\item[$(i)$] $\mu_{n_k}\to \mu$\quad as\;\,$k\to
%\infty$,\;\;as Radon measure on\;\,$M$\,;
%\item[$(ii)$] $V_{\mu_{n_k}}\to V_\mu$\quad as\;\,$k\to
%\infty$,\;\;as Radon measure on\;\,$G_k(M)$\,;
%\item[$(iii)$] $\del V_{\mu_{n_k}}(Y)\to \del V_{\mu}(Y)$\quad as\;\,$k\to
%\infty$,\;\; for any $Y\in C^1_c(TM)$\,;
%\item[$(iv)$] $|\del V_\mu|\le \liminf_{n\to \infty}|\del V_{\mu_{n_k}}|$\;\; as Radon
%measure\,.
%\end{itemize}
%\end{theorem}

\section{Passing measures to
limits}\setcounter{equation}{0}\label{pml} In this section we
address compactness of the family of Radon measures
$\{\mu^\e_t\}_{\e>0}$. At first, recall the following lemma (see Lemma
6.6 in \cite{Ilm2}).
\begin{lemma}\label{lemma14}
Let $\varphi\in C^2_c(M;[0,\infty))$. Then
\[\frac{|\nabla \varphi|^2}{\varphi}\le 2 \max_{\{\varphi>0\}}|\Hess \varphi|\quad \textrm{in}\;\;\{\varphi>0\}\,.\]
\end{lemma}

The next key-fact is known as the {\em semidecreasing property} for the family of measures $\{\mu^\e_t\}_{\e>0}$.
\begin{lemma}\label{lemma7}
Let assumption $(H_0)$ be satisfied. Let $u^\e$ be the solution to
problem \eqref{e1}-\eqref{e1a}. Suppose that \eqref{e121b} and
\eqref{e170} hold true. Let $\varphi\in C^2(M; \re^+)$ with $\suppo
\varphi$ compact. Then, for any $T>0$, the function
\[t \mapsto \mu^\e_t(\varphi)- C(\varphi) t\quad (t\in (0,T))\]
is nonincreasing, for some constant $C=C(\varphi, T)>0$.
\end{lemma}

\noindent{\it Proof\,.\,\,} By \eqref{e1},
\begin{equation}\label{e174}
\begin{split}
\frac{d}{dt}\int_{M} \e \varphi(x) E^\e(x,t)d \mathcal V(x)=
-\int_{M} \e\varphi\left(-\Delta u^\e+\frac
1{\e^2}f(u^\e)\right)^2d \mathcal V(x) \\+\int_{M} \e\langle
\nabla \varphi, \nabla u^\e\rangle \left(-\Delta u^\e+\frac
1{\e^2}f(u^\e) \right)d \mathcal V(x)\,.\hspace{1 cm}
\end{split}
\end{equation}
Hence
\[\frac{d}{dt}\int_{M} \e \varphi(x) E^\e(x,t)d \mathcal V(x)\]
\[\le -\int_{M}\e\varphi\left(-\Delta u^\e+\frac 1{\e^2}f(u^\e)-\frac{\langle \nabla \varphi, \nabla u^\e\rangle}{2\varphi}\right)^2d \mathcal V(x)\]
\[+\int_{M}\e |\nabla u^\e|^2\frac{|\nabla\varphi|^2}{4\varphi}d \mathcal V(x)\le C_1(\varphi)\mu^\e_t(\{\varphi>0\}) \le C\]
for all $t\in (0,T)$, for some positive constant $C=C(\varphi,
T)$ in view of Lemma \ref{lemma14} and the uniform bound \eqref{e170}.
%; here condition $(G^\e_2)$, up to a finite cover, and Lemma
%\ref{lemma14} have been used.
 So, the conclusion follows.\hfill
$\square$

\bigskip

Now, we define the kernel $\phi$ used in \eqref{e60}. In fact, let
$K\subset M$ be a compact subset, $y\in K, s>0$. Let $\hat
\zeta\in C^2([0,\infty))$ such that
\begin{equation}\label{e11}
|\hat \zeta|\le 1\,,\;\; |\hat \zeta'|\le 1\,,\, |\hat \zeta''|\le
1\quad \textrm{in}\;\; [0,\infty)\,,
\end{equation}
\begin{equation}\label{e12}
\hat \zeta = \left\{
\begin{array}{ll}
1 & \textrm{in} \;\;  [0,R_0^2/4)
\\& \\
0  & \textrm{in} \;\; [R_0^2 , \infty)\,,
\end{array}
\right.
\end{equation}
where $R_0:= \inf_{y\in K}\inj(y)\,.$ Define
\[\hat\eta(\rho,t):=[(s-t)]^{-\frac{N-1}2} e^{-\frac{\rho}{4(s-t)}}\quad (\rho\geq 0, 0\le t<s)\,.\]
For any $x\in M$, let
\begin{equation}\label{u52}
\eta(x,t):= \hat \eta(d^2(x), t)\quad (x\in M, 0\leq t<s)\,.
\end{equation}
\[\zeta(x):= \hat \zeta (d^2(x))\quad (x\in M)\,.\]\smallskip
Finally, define
\begin{equation}\label{u52b}
\phi(x,t)\equiv \phi(x,t; y,s):=\eta(x,t)\zeta(x)\quad (x\in M,
0\leq t<s)\,.
\end{equation}

\medskip

In view of \eqref{e170} and the above monotonicity property in Lemma
\ref{lemma7} we can repeat the argument in \cite{Brakke} (see also
\cite{Ilm1}), to show the next compactness result.

\begin{proposition}\label{prop4}
Let assumption $(H_0)$ be satisfied. Let $u^\e$ be the solution to
problem \eqref{e1}-\eqref{e1a}. Suppose that \eqref{e121b} and
\eqref{e170} hold true. Then there are a Radon measure $\mu_t$ on
$M$ and a sequence $\{\e_n\}\subset (0, \infty),\, \e_n\to 0$ as
$n\to\infty$ such that, for every $t>0$,
\begin{equation}\label{e164}
\mu^{\e_n}_t\to \mu_t\quad \textrm{for all}\;\; t\ge
0\;\;\textrm{as}\,\, n\to \infty
\end{equation}
as Radon measure on $M$. Furthermore, for each compact subset
$K\subset M$ we have
\[\int_{M}\phi_{y,s}(x,t) d\mu_t(x)\le \underline C\quad \textrm{for all}\;\; y\in K, 0\le t<s\,, \leqno (G_1)\]
and
\[\mu_t(B_R(x))\le \omega_{N-1} D_0 R^{N-1}\quad \textrm{for all}\;\; x\in K, 0< R<\tilde R, t\ge 0\,. \leqno (G_2)\]
\end{proposition}

\noindent{\it Proof\,.} In view of $(G^\e_2)$ and the above
monotonicity property we can repeat the argument in \cite{Brakke}
(see also \cite{Ilm1}), to show \eqref{e164}. Furthermore, under
hypothesis $(H_0)$, as a consequence of $(G_1^\e)$ and $(G_2^\e)$ from the introduction,
we get $(G_1)$ and $(G_2)$. \hfill $\square$

\section{Clearing-Out Lemma}\setcounter{equation}{0}\label{col}
In this section we will prove the {\it Clearing-Out Lemma}. This result
roughly says that energy concentration occurs only near the
interface region, e.g. $\left\{|u^\e|\leq \alpha \right\}$. In
particular, we show that as $\e \to 0^+$ the solution $u^\e$
converges to either $1$ or $-1$ in a neighborhood of any point
which does not belong to $\overline{\bigcup_{t'\ge 0} \suppo
\mu_{t'}\times\{t'\}}$, where $\mu_t$ is the limit Radon measure
obtained in Section \ref{pml}.

\medskip

For each $y\in M, s>0$ we shall write
\[\phi(x,t;y,s)\equiv \phi_{y,s}(x,t), \quad x\in M, 0\le t<s\,.\]
Observe that because of \eqref{u52b} we clearly have $\phi_{y,s}(x,t)=\phi_{x,s}(y,t)\,$.

\smallskip

\begin{lemma}\label{lemma10} Let assumption $(H_0)$ be
satisfied. Let $u^\e$ be the solution to problem
\eqref{e1}-\eqref{e1a}. Suppose that \eqref{e121b} and
\eqref{e170} hold true. Then

$(i)$ there exists $\kappa_1
>0, \kappa_2>0$ depending on $N, R_0$ and $F$ such that, if for some
$s>t\ge 0$,
\begin{equation}\label{e69}
s-t< \kappa_1,
\end{equation}
and
\begin{equation}\label{e67}
\int_{M} \phi_{y,s}(x,t) d\mu_t(x) <\kappa_2,
\end{equation}
then there exists a neighborhood $V\subset M\times[0,\infty)$ of
$(y,s)$ such that
\begin{equation}\label{e168}
|u^\e|\ge \alpha\quad \textrm{in}\;\;\,V\,
\end{equation}
for all $\e>0$ small enough. As a consequence,
\[(y,s)\not\in \overline{\bigcup_{t'\ge 0} \suppo \mu_{t'}\times\{t'\}}\,.\]

$(ii)$ If $(y,s)\not\in \overline{\bigcup_{t'\ge 0} \suppo
\mu_{t'}\times\{t'\}}$, then there exists a neighborhood $U\subset
M\times[0,\infty)$ of $(y,s)$ such that $\{u^{\e_n}\}$ converges
uniformly in $U$ to either $1$ or $-1$ as $n\to \infty$\,.
\end{lemma}

To prove the Clearing-Out Lemma we need the following technical lemma, which parallels Lemma 3.4 in \cite{Ilm1}
and will be used to prove Lemma \ref{lemma10}. The proof is
standard and is omitted for brevity.

Set
\begin{equation}\label{e169}
\phi_y^r(x)=\phi_{y,s}(x,t)\quad\textrm{whenever}\;\;\,
r^2=2(s-t)\,,
\end{equation}
so
\[\phi_y^r(x):=\zeta(x)\frac 1{ r^{N-1}}e^{-\frac{(d(x,y))^2}{2 r^2}} \, , \qquad x\in M\,.\]

\begin{lemma}\label{lemma9}
Let $\mu$ be a measure satisfying for each $K \subset \subset M$
\[\mu(B_R(y))\le \omega_{N-1} \tilde D_0  R^{N-1}\quad \textrm{for all}\;\; y\in K, 0<R<R_0\]
for some $\tilde D_0=\tilde D_0(K)>0$. Then, for some positive constant
$D=D(R_0,K)$, we have:

\smallskip

\noindent $(i)$ $\int_{M} \phi^r_y d\mu \le D\,;$

\noindent $(ii)$ for any $r>0, 0<R<R_0, y\in M$
$$\int_{M\setminus B_R(y)} \phi^r_y d\mu \le 2^{N-1} D e^{-\frac{3 R^2}{8 r^2}};$$

\noindent $(iii)$ for any $\d>0$ and $\bar r>0$ there exists
$\g_2=\g_2(\d, \bar r)>0$ such that for any $r\in (0, \bar r),
d(y, y_1)\le \g_2 r$ we have
\[\int_{M}\phi^r_{y_1} d\mu \le (1+\d)\int_{M} \phi^r_y d\mu + D\d\,;\]

\noindent $(iv)$ for any $R>0, 0<r<R_0$ with $r\le R$ we have:
\[\int_{M} \phi^r_y d\mu \le \left(\frac{R}{r}\right)^{N-1}\int_{M}\phi^R_y d\mu\,;\]

\noindent $(v)$ for any $\d>0$ there exists $\g_3=\g_3(\d)>0$ such
that for any $r>0, 0<R<R_0$ with $1\le \frac R r\le 1 + \g_3$ and
for any $y\in M$ we have:
\[\int_{M} \phi_y^R d\mu \le (1+\d)\int_{M} \phi^r_y d\mu + D \d\,;\]

\noindent $(vi)$ for any $\d>0$ there exists $\a=\a(\d)>0$ such
that for all $r>0, y\in M$,
\[\int_{M} \phi^{\a(\d)r}_y d\mu \le \frac{\mu(B_r)}{\omega_N[\a(\d)]^{N-1}r^{N-1}}+ \d D\,.\]
\end{lemma}

\smallskip

Now we can prove the {\it Clearing-Out Lemma}.

\noindent {\it Proof of Lemma \ref{lemma10}\,.\,\,} $1.$ In view
of Lemma \ref{lemma9}$(iii)-(v)$ we can a find a neighborhood
$U\subset M\times [0,\infty)$ of $(y,s)$ such that for all $(y',
s')\in U$
\begin{equation}\label{e67a}
0<\inf_{U} (s'-t)\leq s'-t <\kappa_1\,,
\end{equation}
\[\int_{M} \phi_{y',s'}(x,t) d\mu_t(x) <\kappa_2\quad \textrm{for all}\;\; (y',s')\in U\,.\]
We may assume $U\subset\subset M\times (t,\infty)$. By Proposition \ref{prop4} and $(G^\e_1)$,
\begin{equation}\label{e68}
\int_{M} \phi_{y',s'}(x,t) d\mu^{\e_n}_t(x) \le 3\kappa_2\,
\end{equation}
for all $n>n_0$ (for some $n_0\in \ene$) and for all $(y',s')\in
U$.
$4.$  By \cite[Proposition 3.2]{PiPu1}, for any compact
subset $K\subset M$ and for any $\t\in (0,T)$ there exists a
constant $\tilde k>0$ independent of $\e$ such that
\begin{equation}\label{e163}
\|\nabla u^\e(\cdot, t)\|_{L^\infty(K)}\le \frac{\tilde
k}{\e}\quad \textrm{for all}\;\; t\in (\t,T)\,.
\end{equation}
We claim that, for some $C>0$, if $d(x,y')\leq\rho$, then
\begin{equation}\label{e1r}
[F(u^{\e}(y',s'))]^2\leq C \int_M \phi_{y', s'+\frac{\rho^2}2}
(x,s')d\mu^{\e}_{s'}(x),
\end{equation}
where
\[\rho\equiv \rho(y',s'):=\frac{\e}{2 L \tilde k}F(u^{\e}(y',s'))\,,\]
and $L$ is the Lipschitz constant of $F$ in $[-k_0, k_0].$

In fact, whenever $d(x,y')\leq \rho$, we have
\[F(u^{\e}(y',s'))= F(u^{\e}(x,s'))+ F(u^{\e}(y',s'))-F(u^{\e}(x,s'))\]\[\leq F(u^{\e}(x,s'))+ L|u^{\e}(y',s')-u^{\e}(x,s')| \leq F(u^{\e}(x,s'))+
L \tilde k d(x,y') \]
\[\leq F(u^{\e}(x,s'))+ \frac{L\tilde k}{\e}\rho\,.\]
Hence
\[F(u^{\e}(y',s'))\leq 2 F(u^{\e}(x,s'))\,.\]
This easily implies that, for some $C>0,$
\[F(u^{\e}(y',s'))\leq \frac{C}{\rho^N}\int_{B_{\rho}(y')}F(u^{\e}(x,s'))d\mathcal V(x)=\frac C{\rho^{N-1}}\int_{B_{\rho}(y')}
\frac{F(u^{\e}(x,s'))}{\rho}d\mathcal V(x)\,.\] So,
\[[F(u^{\e}(y',s')]^2\leq \frac C{\rho^{N-1}}\int_{B_{\rho}(y')}\frac{F(u^{\e}(x,s'))}{\e}d\mathcal V(x)\leq C \int_M \phi_{y',s'+\frac{\rho^2}2}(x,s')d\mu^{\e}_{s'}(x)\,.\]

By \eqref{e61},
\begin{equation}\label{e69}
\begin{split}
\int_{M} \phi_{y', s'+\frac{\rho^2}2}(x, s') d\mu^{\e}_{s'}(x)\hspace{.7 cm}\\
\le e^{\frac{C_3}2 (\sqrt{s'+\frac{\rho^2}2 - t} -
\frac{\rho}{\sqrt 2} )}\Big[ \int_{M}\phi_{y',
s'+\frac{\rho^2}2}(x,t)d\mu^{\e}_t(x)
+ C_4\Big(s'-t\Big)  + \\ C_5 \left (\sqrt{s'+\frac{\rho^2}2 - t}
- \frac{\rho}{\sqrt 2} \right) \Big]. \hspace{.5 cm}
\end{split}
\end{equation}
Now if $r^2=2(s'-t), R^2=2\left(s'+\frac{\rho^2}2-t\right)$, in view of \eqref{e67a} and $\rho=O(\e)$, then 
$\frac R r\geq 1,\, \frac R r \to 1$ as  $\e\to 0$ uniformly in $U\,.$ So, from Lemma \ref{lemma9}-$(v)$ with $\mu=\mu^\e_t$, for any $\d>0$, for $\e>0$ small enough, we have
\begin{equation}\label{re15}
\int_M \phi_{y', s'+\frac{\rho^2}2}(x,t) d\mu^\e_t \leq (1+\d)\int_M \phi_{y', s'}(x,t)d\mu^\e_t+ D\delta\,.
\end{equation}

In view of \eqref{e1r}, \eqref{e68} and \eqref{re15}, we get
\begin{equation}\label{e3r}
F(u^{\e}(y',s'))\leq e^{C_3\sqrt{\kappa_1}} \big[(1+\delta)\kappa_2+D\d + C_4\kappa_1+C_5\sqrt{\kappa_1}\big]\,.
\end{equation}
If we select $\kappa_1>0, \kappa_2>0$ and $\d>0$ small enough, from
\eqref{e3r} and $(H_0)-(iv)$ it follows that
\begin{equation}\label{e4r}
|u^{\e}| \ge \alpha
\end{equation}
for all $(y',s')\in U$ and for any $\e>0$ sufficiently small.

$4.$ By the same arguments as in \cite[Lemma 3.5]{PiPu1}, there
exists a neighborhood of $(y,s)$,
$V:=B_{R}(y)\times I\subset U$,
such that
\begin{equation}\label{e5r}
1-\e \leq |u^{\e}|\le 1 + \e \quad \textrm{in}\;\; V\,.\end{equation} Let
$s'\in I$. From \eqref{u13} we have that for any $\d'>0$ there
exists $\e_{\d'}>0$ such that for any $\e\in (0, \e_{\d'})$ there holds:
\[\mu^\e_{s'}\big(B_{\frac R2}\big)=\int_{B_\frac R 2}\left[\frac {\e}2|\nabla u^\e|^2+\frac 1{\e} F(u^\e)\right]d \mathcal V(x) \]
\[\le \int_{B_{\frac R2}} \left[ \frac 2{\e} F(u^\e)+\d'\right] d \mathcal
V(x)\le \int_{B_{\frac R2}} \left(\frac 2 {\e} F(1-\e)+ \frac2{\e}
F(1+\e) +\d' \right) d\mathcal V (x)\]\[\le \int_{B_{\frac R2}}
\left[\frac C {\e} \big(\sup_{|s|\leq k_0} |F''(s)|
\e\big)^2+\d'\right] d\mathcal V(x) \le \int_{B_{\frac R 2}}(\bar
C\e+\d') d \mathcal V(x)=\mathcal V(B_{\frac R 2}) (\bar C\e +
\d')\,;\] here use of \eqref{e5r} has been made. Letting $\e\to
0$, and then $\d'\to 0^+$, we get
\[\mu_{s'}\big(B_{R/2}(y)\big)=0\quad \textrm{for}\;\; s'\,\,\textrm{near}\,\, s\,,\]
and
\[u^{\e_n}\to \pm 1\quad \textrm{uniformly in a neighborhood of}\;\; (y,s);\]
thus, $(i)$ has been shown.

\medskip

To prove $(ii)$, let $(y,s)\not\in \overline{\bigcup_{t\ge
0}\suppo \mu_t\times\{t\}}$; hence
\[\int_{M} \phi_{y,s'} d\mu_s'\to 0 \quad \textrm{as}\;\; s' \to s^-,\]
thus,
\[u(y', s')=\lim_{n\to \infty} u^{\e_n}(y', s')=\pm 1\quad \textrm{near}\;\,(y,s)\,.\]
\hfill $\square$

\bigskip
As a simple consequence we have the following result.
\begin{corollary}
We have
\[\suppo \mu\,=\, \overline{\bigcup_{t'\ge 0} \suppo \mu_{t'}\times\{t'\}}\,,\]
where
\[d\mu := d\mu_{t'} dt'\,.\]
\end{corollary}

\noindent{\it Proof\,.\,\,} It is obvious that $\suppo
\mu\,\subseteq\, \overline{\bigcup_{t'\ge 0} \suppo
\mu_{t'}\times\{t'\}}$. Now, let $(y,s)\not\in \suppo \mu$. Then
we can find an open subset $U\subset M\times [0,\infty)$ such that
$(y,s)\in U, \,U\cap \suppo \mu=\emptyset$. Hence $\int_{M}
\phi_{y,s}(x,t)d\mu_t\to 0$ as $t\to s^-\,.$ By Lemma
\ref{lemma10}$(i)$, $(y,s)\not\in \overline{\bigcup_{t'\ge 0}
\suppo \mu_{t'}\times\{t'\}}$; this completes the proof. \hfill
$\square$

\bigskip

Another consequence of the {\em Clearing-Out Lemma} is that one can control the size of the set where the energy is concentrating.

\begin{corollary}\label{corhauss}
Let $U\subseteq M$ be an open subset. Then there exists $C_5>0$
and $C_6>0$ such that if $(\suppo \mu)_t:=\suppo \mu \cap \big(M\times\{t\}\big)\,$ ,
then
\medskip

\noindent $(i)$ $\mathcal H^{N-1}\big((\suppo \mu)_t \cap U\big)
\le C_5\liminf_{s\to t^-} \mu_s (U)$ for every $t>0,$

\noindent $(ii)$ $\mathcal H^{N-1}\big((\suppo \mu)_t \cap
B_R\big)\le C_6 R^{N-1}$ for every $0<R<R_0$ and $t\ge 0$\,.
\end{corollary}

\noindent{\it Proof\,.\,\,} Clearly, $(i)$ follows, if we show
that $$\mathcal H^{N-1}\big((\suppo \mu)_t \cap K\big) \le
C\liminf_{s\to t^-} \mu_s (U)$$ for every $t>0,$ for every compact
subset $K\subset U$.

Let $(y,t)\in (\suppo \mu)_t\cap K;$ take any $\d>0$ and
$\a=\a(\d)>0$ given by Lemma \ref{lemma9}$(vi)$. For each
$0<r<\sqrt{\kappa_1}$, by Lemma \ref{lemma10}
\begin{equation}\label{e77}
\kappa_2 \le \int_{M} \phi^{\a r}_y d\mu_{t-\frac{\a^2r^2}{2}}\,.
\end{equation}
which combined with Lemma \ref{lemma9}$-(vi)$ yields
\[\kappa_2 \le \frac{\mu_{t-\frac{\a^2r^2}{2}}(B_r)}{\omega_N r^{N-1}[\a(\d)]^{N-1}}\, + D\d\,.\]
The previous inequality with $\delta=\frac{\kappa_2}{2 D}$ gives
\begin{equation}\label{e78}
\kappa_2\le \frac{2 \mu_{t-\frac{\a^2r^2}{2}}(B_r)}{\omega_N
r^{N-1}[\a(\d)]^{N-1}}\,.
\end{equation}
Let $r>0$ such that $dist(K, \pa U)>r$ and consider the covering
of $\suppo \mu_t\cap K$
\[\mathcal B = \{B_r(x)\,:\,x\in \,(\suppo \mu)_t\,\}\,.\]
By the Besicovitch Covering Theorem, which can be applied for
compact subset of Riemannian manifolds, (see, $e.g.$, Theorem
1.1.4 and Example 1.15 $(c)$ in \cite{Hein}), there are finitely
many countable subcollections $\mathcal B_1, \ldots, \mathcal
B_{\bar l}$ of $\mathcal B$ such that each $\mathcal B_i$ is made
of disjoint balls, and
\[(\suppo \mu)_t \cap K \subseteq \bigcup_{l=1}^{\bar l}\bigcup_{B_r(y_j)\in \mathcal B_l}\, B_r(y_j)\,.\]
We have, for some $\tilde C>0$, an estimate for the pre-Hausdorff
measures
\[\mathcal H_r^{N-1}\Big((\suppo \mu)_t \cap K\Big)\le \tilde C \sum_{l=1}^{\bar l}\sum_{B_r(y_j)\in \mathcal
B_i}\omega_N r^{N-1} \]
\[
\le \tilde C\sum_{l=1}^{\bar
l}\frac2{\a^{N-1}\kappa_2}\sum_{B_r(y_j)\in \mathcal B_i}
\mu_{t-\a^2r^2/2}\big(B_r(y_j)\big)
\]
\[\le \tilde C\sum_{l=1}^{\bar l}\frac
2{\a^{N-1}\kappa_2}\mu_{t-\a^2r^2/2}\Big(\big\{x : dist (x, K)\le
r\big\} \Big)\le \frac{2 \tilde C\bar
l}{\a^{N-1}\kappa_2}\mu_{t-\a^2r^2/2}(U).
\]
Sending $r\to 0^+$, we obtain
\[\mathcal H^{N-1}\Big((\suppo \mu)_t\cap K\Big)\le C \liminf_{s\to t^-}\mu_s (U)\,,\]
so, $(i)$ has been proven. Furthermore, $(ii)$ follows by $(i)$
and $(G_2)$. \hfill $\square$

\begin{remark}\label{emptysp}
Notice that the {\it Clearing-out
Lemma} could also be proved analogously to \cite{Ilm1}. In that
case, we need the so-called {\it Empty Spot Lemma}, which could be
deduced in the present situation, too. Indeed, it is mainly based
on a result given in \cite{Chen} in Euclidean space, concerning
propagation of interfaces, that could be easily shown also in a general Riemannian manifold $M$.
To do this the key role is played by an important property of the the signed distance $\tilde d(x,t)$ from $\pa \Sigma^t$,
$\Sigma^t$ being a family of sets evolving by mean curvature flow,
starting from a sphere $\pa B_R(x_0)$, with $R>0$ small enough. Indeed, as in \cite[Section 7]{AAM}, in a tubular neighborhood 
of $\Sigma^t$, one has:
\begin{equation}\label{e60a}
\big| \pa_t \tilde d(x,t) - \Delta \tilde d(x,t)\big| \le
\tilde C |\tilde d(x,t)|\,.
\end{equation} This is all what is needed to conclude. 
\end{remark}
\medskip

\section{Density lower bound}\setcounter{equation}{0} \label{dlb}
The result of this section, roughly speaking, shows that a suitable
$(N-1)-$density of $\mu$ ({\em forward density} in the terminology of \cite{Ilm1}) is bounded below on the support of the
measure. More precisely, we show that an explicit lower bound
holds $\mathcal H^{N-1}-$a.e. on each time-slice of $\suppo
\mu\,.$

\smallskip

In the sequel, we take $\kappa_2$ as in Lemma \ref{lemma10}.
Define
%Define
%\[Z^0:=\left\{(x,t)\in \suppo\mu\,:\, \limsup_{r\to 0^+}\int_{M} \phi^r_x d\mu_t\, < \, \kappa_2\right\}\,,\]
%\[Z^0_t:= Z^0\cap \Big(M\times\{t\} \Big)\,.\]

%\begin{lemma}\label{lemma11}
%For any $\s>0$, $\mathcal H^{N-2+\s}(Z^0_t)=0$ for a.e. $t\ge
%0$\,.
%\end{lemma}
\[Z^-:=\left\{(x,t)\in \suppo \mu\,:\,\limsup_{s\to t^+} \int_{M} \phi_{y,s} d\mu_s(y)\,<\,\kappa_2 \right\}\,,\]
\[Z^-_t:=Z^-\cap \Big(M\times\{t\} \Big)\,. \]

\begin{lemma}\label{lemma12}
For any $\tilde \s>0$, $\mathcal H^{N-2+\tilde\s}(Z^-_t)=0$ for
a.e. $t\ge 0$\,.
\end{lemma}

\noindent{\it Proof\,.\,\,} $1.$\, It is direct to see that, for
each $\bar \t>0,$
\[Z^- = \bigcup_{0<\t<\bar \t, \kappa_3<\kappa_2}Z^{\kappa_3, \t}\,,\]
where
\[Z^{\kappa_3, \t}:=\left\{(x,t)\in \suppo\mu\,:\,\int_{M} \phi_{y,s}d\mu_s(y) \le \kappa_3\,\,\textrm{for all}\,\,s\in(t, t+\t)\right\}.\]
Hence, the thesis will follow, if we prove that $\mathcal
H^{N-2+\tilde\s}(Z_t^{\kappa_3, \t})=0$ for each
$\kappa_3<\kappa_2, 0<\t<\kappa_1$.

\smallskip
\medskip

{\it Claim.} Let $\d>0, s\in [t, t+\t], \g_2=\g_2(\d, \sqrt{2\t})$
be the constant given by Lemma \ref{lemma9}$(iii)$,
$r:=\sqrt{2(s-t)}, t':=s+\frac{r^2}2=t+r^2.$ If
$$t'-t\le 2\t\quad \textrm{and}\;\; d(x,x')\le \g_2 r\,,$$
then
$$\textrm{either}\;\,(x,t)\not \in Z^{\kappa_3, \t}\;\; \textrm{or}\;\; (x',t')\not\in Z^{\kappa_3, \t}\,.$$
Indeed, more precisely, we are going to show that if $(x,t)\in
Z^{\kappa_3, \t}$, then
\begin{equation}\label{e167}
\mathcal P_{2\t}^{(x,t)}\cap Z^{\kappa_3, \t}=\{(x,t)\},
\end{equation}
where
\[\mathcal P_{2\t}^{(x,t)}:=\left\{(x,t)\in M\times [0,\infty)\,:\,\frac{d^2(x',x)}{\g_2}\le t'-t\le 2\t\,\right\}\,.\]

\medskip

In fact, let $(x,t)\in Z^{\kappa_3, \t}.$ By Lemma
\ref{lemma9}$(iii)$, for any $x'\in B_{\g_2 r}(x)$, we have:
\[\int_{M} \phi_{x', s+r^2/2}(y,s) d\mu_s(y) = \int_{M} \phi^r_{x'}(y)d\mu_s(y)\]
\[\le (1+\d)\int_{M} \phi^r_x(y) d\mu_s(y) + D \d=(1+\d)\int_{M} \phi_{y,s}(x,t)d\mu_s(y) + D\d\]
\[\le (1+\d)\kappa_3 + D\d<\kappa_2,\]
for $\d>0$ sufficiently small. By Lemma \ref{lemma10}, $(x',
t')\not \in \overline{\bigcup_{\xi\ge 0} \suppo \mu_{\xi}}$; thus,
$(x',t')\not\in Z^{\kappa_3, \t}.$ This proves the claim.

\smallskip

$2.$\, For every $ x_0\in M, t_0>0$, define
\[Z' :=Z^{\kappa_3, \t}\cap \big( B_1(x_0)\times [t_0-\t, t_0+\t]\big)\,.\]
Since $Z^{\kappa_3, \t}$ is a countable union of such $Z'$, the
thesis follows, if we show that $\mathcal
H^{N-2+\tilde\s}(Z'_t)=0$ for a.e. $t\ge 0$, where $Z^\prime_t:= Z'\cap
\big(M \times\{t\} \big)\,.$

Observe that the set $Z'\cap\big(\{x\}\times \re \big)$ contains
at most one point, for $\mathcal P_{2\t}^{(x,t)}$ is higher than
$Z'$ when $(x,t)\in Z',$ in view of \eqref{e167}.

\smallskip

Let $\pi$ be the (nearest point) projection from $M\times [0,\infty)$
onto $M\times\{0\}$; so $\pi (Z')\subset B_1(x)$. Let $\d_1>0$.
There exist sequences $\{x_i\}_{i\in\ene}\subset \pi (Z')$ and
$\{r_i\}_{i\in\ene}\subset (0,\d_1)$ with
\begin{equation}\label{e78}
\sum_{i=1}^\infty \omega_N r_i^N \le 2 \mathcal
V\big(B_1(x_0)\big)\,,
\end{equation}
such that
\[\bigcup_{i=1}^\infty B_{r_i}(x_i)\supseteq \pi(Z')\,.\]

\smallskip

In view of Step $1.$,
\[Z'\subseteq B_{r_i}(x_i)\times [t_i -r_i^2/\g^2, t_i+ r_i^2/\g^2]\,,\]
where $(x_i, t_i):=\pi^{-1}(x_i).$ Thus, for some $\tilde C>0,$ we
have:
\[\int_{t_0-\t}^{t_0+\t} \mathcal H^{N-2+\tilde\s}_{\d_1}(Z'_t) dt\le \int_{t_0-\t}^{t_0+\t}\sum_{\{i\in \ene :t\in [t_i-r_i^2/\g^2, t_i+r_i^2/\g]\}}
\tilde C r_i^{N-2+\tilde\s}dt\]
\[=\tilde C \sum_{i=1}^\infty \int_{t_i-r_i^2/\g^2}^{t_i + r_i^2/\g^2} r_i^{N-2+\tilde\s}=\tilde C \sum_{i=1}^\infty \frac{2}{\g^2}r_i^{N+\tilde\s}\]
\[\le 2\tilde C D\d_1^{\tilde\s} \mathcal V\big(B_1(x_0)\big)\,.\]
Sending $\d_1\to 0^+,$ by the monotone convergence theorem
\[\int_{t_0-\t}^{t_0+\t} \mathcal H^{N-2+\tilde\s}(Z'_t) dt = 0\,.\]
This implies the result. \hfill $\square$

\section{Vanishing of the limit
discrepancy}\setcounter{equation}{0}\label{vld} The purpose of
this section is to show that the discrepancy Radon measure
vanishes as $\e\to 0^+,$ up to subsequences. More precisely,
define
\[ d\xi^\e = d\xi^\e_t dt, \;\;\; d\mu^\e:= d\mu^\e_t dt.\]
Since $|\xi^\e|\le \mu^\e$, by Proposition \ref{prop4} we can
assume that there exist a Radon measure $\xi$ on $M\times
[0,\infty)$, and a subsequence of $\{\e_n\}$, which will be still
denoted by $\{\e_n\}$, such that
\[\xi^{\e_n} \to \xi,\;\; \mu^{\e_n}\to \mu:= d\mu_t dt\quad \textrm{as}\;\; n\to \infty\]
as Radon measure on $M\times [0,\infty)$. By \eqref{u13}, $\xi \le
0$. Indeed, we are going to show the following result.

\begin{proposition}\label{prop6}
There holds\, $\xi = 0$\,.
\end{proposition}

In order to prove Proposition \ref{prop6} we use the following
lemma (see \cite{PiPu1}).

\begin{lemma}\label{lemma5}
Let assumption $(H_0)$ be satisfied. Let $u^\e$ be the solution to
problem \eqref{e1}-\eqref{e1a}. Suppose that \eqref{e121b} and
\eqref{e170} hold true. Let $K\subset M$ be a compact subset,
$y\in K, s>0.$ Let $\phi:=\eta\zeta$ with $\zeta$ and $\eta$ as in
Section \ref{pml}. Then for every $\e>0$
\begin{equation}\label{e47}
\begin{split}
\frac{d}{dt}\int_{M} \phi(x,t)d\mu^\e_t(x)\,\le\, \frac 1{2(s-t)}
\int_{M} \phi\e\big\{|\nabla u^\e|^2-E^\e \big\}d \mathcal V(x)
\\+\frac {C_3}{(s-t)^{1/2}}\int_{M} \phi
d\mu^\e_t(x)+ C_4\quad \textrm{for
all}\;\;0<t<s\,.\hspace{.3 cm}
\end{split}
\end{equation}
for some positive constants $C_3$ and $C_4$ only dependending on
$K$.
\end{lemma}

\noindent{\it Proof of Proposition \ref{prop6}\,.\,\,} {\it Step}
$1.$\, Given $x_0\in M$ and $K=\overline{B_{R_0}(x_0)}$, in view
of $(G^{\e}_1)$, integrating \eqref{e47} we get
\[- \int_0^{s-\tilde\s} \int_{M} \frac 1{2(s-t)} \phi_{y,s}(x,t) d\xi^{\e}_t(x) dt \le \int_{M} \phi_{y,s}(x,0) d\mu^{\e}_0(x)\]
\[+ C_3 \underline C \sqrt s + C_4  s \]
for every $(y,s)\in M\times(0,\infty), 0<\tilde\s<s\,.$ Letting
$\e\to 0$, using Lemma \ref{lemma9}$(i)$ and $(G_1)$, since
$\xi\le 0,$ we obtain a uniform bound, for $s$ in a compact set,
namely
\[\int_0^{s-\tilde\s} \int_{M} \frac 1{2(s-t)} \phi_{y,s}(x,t) d|\xi|(x,t) \le \int_{M} \phi_{y,s}(x,0) d\mu_0(x)\]
\[+ C_3 \underline C \sqrt s + C_4 s\le \underline C + C_3 \underline C \sqrt s + C_4  s=: \bar C\,.\]
\smallskip

Now, take $T>0, 0<R<R_0$. Integrating over $B_R(x_0)\times (0,
T+1)$ we have:
\[\int_{\tilde\s}^{T+1}\int_{B_R(x_0)}\int_0^{s-\tilde\s}\int_{M} \frac 1{2(s-t)}
\phi_{y,s}(x,t) d|\xi|(x,t)d\mu_s(y) ds\]\[ \le\bar C
\omega_{N-1}D_0(T+1)R^{N-1}\,.\] Hence, by Tonelli theorem,
\[\int_0^{T+1-\tilde\s}\int_{M}\int_{t+\tilde\s}^{T+1}\frac 1{2(s-t)}\int_{B_R(x_0)} \phi_{y,s}(x,t)
d\mu_s(y) ds d|\xi|(x,t)\]\[\le\bar C
\omega_{N-1}D_0(T+1)R^{N-1}\,.
\] Sending $\tilde\s\to 0$, by the monotone convergence theorem we
have
\[ \int_0^{T+1}\int_{M}\int_{t}^{T+1}\frac 1{2(s-t)}\int_{B_R(x_0)} \phi_{y,s}(x,t) d\mu_s(y) ds d|\xi|(x,t)\]\[\le\bar C \omega_{N-1}D_0(T+1)R^{N-1}\,.\]
So,
\[\int_t^{t+1} \frac 1{2(s-t)}\int_{B_R(x_0)}\phi_{y,s}(x,t)d\mu_s(y) ds \le C(x,t)<\infty\]
for $|\xi|-$a.e. $(x,t)\in M\times (0,T)\,.$

\medskip
\smallskip

{\it Step} $2.$ For any $x\in B_{R/2}(x_0), s>t>0$ we have:
\[\int_M \phi_{y,s}(x,t)d\mu_s(y) = \int_{B_R(x_0)}\phi_{y,s}(x,t) d\mu_s(y) + \int_{M\setminus B_R(x_0)}\phi_y^{\sqrt{2(s-t)}}(x)d\mu_s(y)\]
\[\le \int_{B_R(x_0)}\phi_{y,s}(x,t) d\mu_s(y) + \int_{M\setminus
B_{R/2}(x_0)}\phi_y^{\sqrt{2(s-t)}}(x) d\mu_s(y)\]\[\le
\int_{B_R(x)}\phi_{y,s}(x,t) d\mu_s(y) + 2^{N-1}e^{-\frac 3 8
\frac{(R/2)^2}{2(s-t)}}D;\] here the fact that $B_{R/2}(x)\subset
B_R(x_0)$ for any $x\in B_{R/2}(x_0)$, and Lemma \ref{lemma9}-(ii)
have been used. Thus for $|\xi|-$ a.e. $(x,t)\in
B_{R/2}(x_0)\times [0,T]$,
\begin{equation}\label{e79}
\begin{split}
\int_t^{t+1}\frac{1}{2(s-t)}\int_{M} \phi_{y,s}(x,t)d\mu_s(y) ds
\le \hspace{.6 cm}\\ \le C(x,t)+ \int_t^{t+1}\frac
1{2(s-t)}2^{N-1}e^{-\frac 3{32}\frac {R^2}{2(s-t)}} D<\infty\,.
\end{split}
\end{equation}
Since $T>0$ and $x_0\in M$ were arbitrary, \eqref{e79} holds for
$|\xi|-$a.e. $(x,t)\in M\times [0,\infty)$.

\smallskip
\medskip

{\it Step} $3.$ Now, take $(x,t)\in M\times [0,T]$ such that
\eqref{e79} holds true. We shall prove that
\[\lim_{s\to t^+} \int_{M} \phi_{y,s}(x,t) d\mu_s(y) \,=\,0\,.\]
Let $\b(s):=\log (s-t),$ and
\[h(s):=\int_{M} \phi_{y,s}(x,t) d\mu_s(y)\,.\]
Hence, \eqref{e79} implies
\begin{equation}\label{e80}
\int_{-\infty}^0 h(t+e^\b) d\b < \infty\,.
\end{equation}
We are going to show that \eqref{e80} yields $h(t+e^\b)\to 0$ as
$\b\to-\infty$\,.
\medskip
\smallskip

{\it Step} $4.$\, Let $\g\in (0,1]$ to be specified later. By
\eqref{e80}, there exists a sequence $\{\b_n\}\subset (-\infty,0)$
such that
\begin{equation}\label{e81}
\lim_{n\to \infty} \b_n=-\infty,\;\,0<\b_n - \b_{n+1}\le \g,\;\,
h(t+e^{\b_n}) \le \g\,.
\end{equation}
Let $\b\in (-\infty, \b_1]$ and assume $\b\in [\b_n, \b_{n-1})$
for some $n\in \ene$. Since $\b_n\le \b$, from
\eqref{e61} we get
\begin{equation}\label{e82}
\begin{split}
h(t+ e^{\b})=\int_{M} \phi_{y, t+e^\b}(x,t)
d\mu_{t+e^\b}(y)\hspace{2 cm}\\=\int_{M} \phi_{x, t+e^{2\b}}(y,
t+e^\b)d\mu_{t+e^\b}(y)\hspace{2.2 cm}
\\ \le
e^{\frac{C_3}2 \big(\sqrt{e^{2\b}-e^{\b_n}} -  \sqrt{e^{2\b}-e^{\b}}\big)} \Big[
\int_{M}\phi_{x, t+e^{2\b}}(y,
t+e^{\b_n})d\mu_{t+e^{\b_n}}(y)\\ +C_4 (e^\b-e^{\b_n})   +   C_5 \big(\sqrt{e^{2\b}-e^{\b_n}} -  \sqrt{e^{2\b}-e^{\b}}\big)       \Big]\\
= e^{\frac{C_3}2 \big(\sqrt{e^{2\b}-e^{\b_n}} -  \sqrt{e^{2\b}-e^{\b}}\big)}  \Big[
\int_{M}\phi_x^R(x, t+e^{\b_n})
d\mu_{t+e^{\b_n}}(y) \\ +C_4 (e^\b-e^{\b_n})   +   C_5 \big(\sqrt{e^{2\b}-e^{\b_n}} -  \sqrt{e^{2\b}-e^{\b}}\big)\Big]\,,\hspace{.5
cm}
\end{split}
\end{equation}
where $R:=\sqrt{2(2e^\b-e^{\b_n})}\,.$ Furthermore, in view of
\eqref{e81}, we have
\begin{equation}\label{e83}
\g \ge h(t+ e^{\b_n})=\int_{M} \phi_{y, t+ e^{\b_n}}(x,t) d\mu_{t+
e^{\b_n}}(y)=\int_{M}\phi_x^r(y) d\mu_{t+ e^\b_n}(y)\,,
\end{equation}
where $r:=\sqrt{2 e^{\b_n}}\,.$ Note that
\begin{equation}\label{e84}
1 \le \frac R r=\sqrt{2 e^{\b-\b_n}-1}\le 1 + \tilde C\g\,.
\end{equation}

\smallskip
\medskip

{\it Step} $5.$\, Let $\d>0$ and set $\g=\min\{\d, \g_3(\d)/\tilde
C\},$ where $\g_3$ is given by Lemma \ref{lemma9}$(iv)$. By
\eqref{e82}-\eqref{e83} and Lemma \ref{lemma9}$(v)$,
\[ h(t+ e^\b)\le
e^{\frac{C_3}2 \big(\sqrt{e^{2\b}-e^{\b_n}} -  \sqrt{e^{2\b}-e^{\b}}\big)}     \Big[
\int_{M}\phi_x^R
d\mu_{t+e^{\b_n}}(y)\]\[ +C_4 (e^\b-e^{\b_n})   +   C_5 \big(\sqrt{e^{2\b}-e^{\b_n}} -  \sqrt{e^{2\b}-e^{\b}}\big)\Big]\]
\[\le   e^{\frac{C_3}2 \big(\sqrt{e^{2\b}-e^{\b_n}} -  \sqrt{e^{2\b}-e^{\b}}\big)}    \Big[
(1+\d)\int_{M}\phi_x^r d\mu_{t+e^{\b_n}}(y)+D \d\]\[+   C_4
(e^\b-e^{\b_n})   +   C_5 \big(\sqrt{e^{2\b}-e^{\b_n}} -
\sqrt{e^{2\b}-e^{\b}}\big)    \Big]\]\[\le  e^{\frac{C_3}2
\big(\sqrt{e^{2\b}-e^{\b_n}} -  \sqrt{e^{2\b}-e^{\b}}\big)} \Big[
(1+\d)\g+ D\d\]\[+  C_4 (e^\b-e^{\b_n})   +   C_5
\big(\sqrt{e^{2\b}-e^{\b_n}} -  \sqrt{e^{2\b}-e^{\b}}\big) \big]\]
for all $\b\in [\b_{n}, \b_{n-1})$\,. So, letting $\d\to 0^+$ (and
thus $\g\to 0^+$), as $\b\to-\infty$ (hence $\b_n\to-\infty$) we
obtain
\begin{equation}\label{e85}
\lim_{s\to t^+} h(s)=0 \quad \textrm{for}\;\;
|\xi|-\textrm{a.e.}\,\, (x,t)\in M\times (0,T)\,.
\end{equation}

\smallskip
\medskip

{\it Step} $6.$ By Lemma \ref{lemma12}, for any $\tilde\s>0$,
\begin{equation}\label{e86}
\limsup_{s\to t^+} h(s) \ge \kappa_2 >0\quad d\big(\mathcal
H^{N-2+\tilde\s}\lfloor \suppo \mu_t\,\big) dt-\textrm{a.e.}\,
(x,t)\in M\times (0,\infty)\,.
\end{equation}
On the other hand, in view of $(G_2)$, for all $\bar x\in M,
0<R<R_0$, in $B_R(\bar x)\times [0,T]$ for any $\tilde\s \in
(0,1)$ there holds:
\begin{equation}\label{e87}
d|\xi| \le d\mu= d\mu_t dt << d\big(\mathcal
H^{N-2+\tilde\s}\lfloor \suppo \mu_t\,\big) dt\,.
\end{equation}
By \eqref{e85}-\eqref{e87},
\[0<\kappa_2\le \limsup_{s\to t^+} h(s) = 0\quad |\xi|-\,\textrm{a.e.}\, (x,t)\in M\times (0,T).\]
This implies $\xi\lfloor_{B_R(\bar x)\times (0,T)}=0.$ Since $\bar
x>0$ and $T>0$ were arbitrary, the conclusion follows. \hfill
$\square$

\section{Brakke's inequality for the limit
measure}\setcounter{equation}{0}\label{bilm} In this section we
establish the main result of the present paper. In fact, we prove
that the limit measure $\mu_t$ evolves by mean curvature flow, in
the sense of {\it Brakke}. To state this result precisely, we need
some notations.

Recall that the {\it upper derivative} of a function $\psi:\re\to
\re$ at $x_0\in \re$ is given by
$$\ud_{x_0} \psi:=\limsup_{x\to x_0} \frac{\psi(x)- \psi(x_0)}{x-x_0}\,.$$

Let $\phi\in C^2_c(M ; \re^+)$; for any $\mu\in \mathcal M(M)$
define:
\[\mathfrak{B}(\mu, \phi)\equiv -\infty\]
whenever either of the following holds:
\begin{itemize}
\item[$(i)$]\,\; $\mu\lfloor\{\phi>0\}\not\in\mathcal M_k(M)$;
\item[$(ii)$]\, $|\d\mu|\lfloor\{\phi>0\}$ is not absolutely
continuous with respect to $\mu\lfloor\{\phi>0\}$\,;
\item[$(iii)$] $\int_{M} \phi H^2 d\mu =\infty\,.$
\end{itemize}
Otherwise,
\begin{equation}\label{e171}
 \mathfrak{B}(\mu,\phi):=\int_{M}\big\{ -\phi H^2 + \langle \nabla \phi,
T_x\mu^{\perp}(\,\overrightarrow{H}\,)\rangle\,\big\}\,d\mu\,.
\end{equation}

Define
\[\mathfrak{B}^{\e}(u^\e, \phi):= -\e\int_{M} \left\{\phi\left(-\Delta u^{\e} +\frac{f(u^\e)}{\e^2}\right)^2+
\langle \nabla \phi, \nabla u^\e\rangle\left(-\Delta u^\e +\frac
{f(u^\e)}{\e^2}\right)\right\}d \mathcal V\,\] and observe that, in
view of \eqref{e174}, for any $t>0$ we have
\begin{equation}\label{e175}
\frac{d}{dt}\mu^\e_t(\phi)=\mathfrak{B}^{\e}(u^\e, \phi)\,.
\end{equation}

Our purpose is pass to the limit as $\e \to 0^+$ in \eqref{e175}.
Indeed, to do it appropriately we shall use suitable {\it
varifolds} (see Subsection \ref{varif}) associated to $u^\e$, and
results proved in Sections \ref{pml}, \ref{dlb}, \ref{vld}. Our
main result is as follows.

\begin{theorem}\label{thmdb}
Let assumptions $(H_0)$ be satisfied. Let $u^\e$ be the solution
to problem \eqref{e1}-\eqref{e1a}. Suppose that \eqref{e121b} and
\eqref{e170} hold true. Then the family of Radon measures
$\{\mu_t\}_{t\ge 0}$ from Proposition \ref{prop4} are
$(N-1)-$rectifiable for a.e. $t>0$ and satisfies the \em{Brakke's inequality}:
\begin{equation}\label{e90}
\ud_t \mu_t(\phi)\, \le\, \mathfrak{B} (\mu_t, \phi)\,
\end{equation}
for every $ \phi\in C^2_c(M ; \re^+)\, $ and for every $t>0$,
where
\[\mu_t(\phi)\equiv \int_{M} \phi(x) d\mu_t(x)\,.\]

\end{theorem}

Before going into the proof of the main result, let us introduce varifolds that will be used in the sequel.

\smallskip

Let $t_0\ge 0; \{t_n\}_{n\in\ene}\subset [0,\infty),\; t_n\to
t_0.$ Let $u^\e$ be a family of equibounded solutions to problem
\eqref{e1}-\eqref{e1a}. Let $\{\e_n\}\subset (0,1),\, \e_n\to 0$;
consider the sequence of functions
\[\{u^n\}\equiv \{u^{\e_n}(\cdot, t_n)\}\,.\]
Define
\[d\mu^n:=\left(\frac{\e_n}2 |\nabla u^n|^2+\frac{F(u^n)}{\e_n}\right) d \mathcal V\,;\]
\[d\xi^n:=\left(\frac{\e_n}{2}|\nabla u^n|^2 -\frac{F(u^n)}{\e_n}\right)d \mathcal V\,.\]
By standard results in unique continuation for parabolic equations, for each $t>0$ and $n\in \ene$,
\begin{equation}\label{e95}
 \mathcal V\big(\{\nabla u^{\e_n}(\cdot,t)=0\}\big)\,=\,0\,.
\end{equation}
So, for all $n\in \ene$, we can define the $(N-1)-$ varifold $V^n$
by
\[\int_{M} \psi (x, S) d V^n(x,S)\,=\,\int_{M} \psi (x, \nabla u^n(x)^{\perp})d\mu^n(x)\]
for any\;\;$\psi\in C^0_c(G_{N-1}(M); \re)$. Note that
$\|V\|=\mu^{n}\,.$ Define
\begin{equation}\label{e172}
\begin{split}
\mathfrak{B}^{n}(u^n, \phi):= -\e_n\int_{M} \Big\{\phi\Big(-\Delta
u^n +\frac{f(u^n)}{\e_n^2}\Big)^2  \\+ \langle \nabla \phi, \nabla
u^n\rangle\Big(-\Delta u^n +\frac {f(u^n)}{\e_n^2}\Big)\Big\}d
\mathcal V\,.\hspace{.5 cm}
\end{split}
\end{equation}

The next proposition will have a key role in the proof of Theorem
\ref{thmdb}
\begin{proposition}\label{prop7}
Let $\phi\in C^2_c(M; \re^+)$. Assume that
\begin{itemize}
\item[$(i)$] there exists $\mu$ such that $\mu^n\to \mu\;$ as $n\to \infty$, as Radon measure on $M$\,;
\item [$(ii)$] $|\xi^n|\big(\{\phi>0\}\big)\to 0$ as
$n\to\infty$\,;
\item[$(iii)$] there exists a constant $\check{C}>0$ such that $\mathfrak B^{n}(u^n, \phi)\ge -
\check C$ for all $n\in \ene$\,;
\item [$(iv)$] $\mathcal H^{N-1}\big(\suppo \mu
\,\cap\,\{\phi>0\}\big)<\infty$\,.
\end{itemize}
Then
\begin{itemize}
\item[$(a)$] $\mu\lfloor\{\phi>0\}$ is $(N-1)-$rectifiable;
\item[$(b)$] there exists $V\in \mathcal R\mathbb V_{N-1}(M)$
such that $V^n\lfloor\{\phi>0\}\to V$ as $n\to \infty,$ and
$\|V\|=\mu\lfloor\{\phi>0\}\,;$
\item[$(c)$] for all $Y\in C^1_c\big(\{\phi>0\}; TM\big),$
\[\del V(Y)= - \lim_{n\to \infty} \int_{M} \e_n \langle Y, \nabla u^n\rangle\left(-\Delta u^n+\frac{f(u^n)}{\e_n^2}\right)d \mathcal V\,;\]
\item[$(d)$] $\mathfrak B(\mu, \phi)\ge \limsup_{n\to \infty} \mathfrak
B^{n}(u^n, \phi)\,.$
\end{itemize}
\end{proposition}

\subsection{Proof of Proposition \ref{prop7}}
The proof of Proposition \ref{prop7}, to which this Subsection is
devoted, needs some preliminary results. To begin with, the next
density lemma will be used (see Lemma 7.4 in \cite{Ilm2}).

\begin{lemma}\label{lemma15}
Let $\mu\in \mathcal M_k(M)$. For any vector field $Z\in
L^2_{\mu}(TM)$ and any $\delta>0$, there exists a vector field
$Y\in C^1_c(TM)$ such that
\[\|Z-Y\|_{L^2_\mu(TM)}\le \delta\,.\]
\end{lemma}

Furthermore, we make use of the following lemma.
\begin{lemma}\label{lemma13}
Let $\phi\in C^2_c(M; \re^+),\, \mu$ be a Radon measure,
$\displaystyle{C_1(\phi):=\sup_{M} |\Hess \phi|\,.}$ Let $u^\e$ be
solution to equation \eqref{e1}\,. Then
\begin{equation}\label{e91}
\int_{M} \langle\nabla \phi, T_x \mu^{\perp}(\,
\overrightarrow{H}\,)\rangle d\mu \le \frac 1 2\int_{M} \phi H^2
d\mu + C_1(\phi)\mu (\{\phi>0\})\,,
\end{equation}
\begin{equation}\label{e92}
\int_{M} \phi H^2 d\mu \le - 2 \mathfrak B (\phi,\mu) + 2
C_1(\phi)\mu(\{\phi>0\})\,,
\end{equation}
when these are defined; similarly
\begin{equation}\label{e93}
\begin{split}
\int_{M} \e \langle\nabla \phi, \nabla u\rangle\left(-\Delta u +
\frac 1{\e^2} f(u^\e) \right)d \mathcal V(x)  \hspace{.9 cm}
\\
\le \frac 1 2\int_{M} \e\phi\left(-\Delta u^\e +\frac 1{\e^2}
f(u^\e)\right)^2 d \mathcal V(x) + 2
C_1(\phi)\int_{\{\phi>0\}}\frac{\e}2|\nabla u^\e|^2 d \mathcal
V(x)\,;
\end{split}
\end{equation}
\begin{equation}\label{e94}
\begin{split}
\int_{M} \e \phi\left(-\Delta u^\e +\frac 1{\e^2} f(u^\e)
\right)^2 d \mathcal V(x)  \hspace{.4 cm} \\ \le - 2 \mathfrak
B^\e(u^\e, \phi) + 4 C_1(\phi)\int_{\{\phi>0\}}\frac{\e}2|\nabla
u^\e|^2 d \mathcal V(x)\,.
\end{split}
\end{equation}
\end{lemma}

\noindent{\it Proof\,.\,\,} Inequality \eqref{e91} follows from
the Cauchy-Schwartz inequality, since $C_1(\phi)\ge
\sup_{M}\frac{|\nabla \phi|^2}{\phi}$, in view of Lemma
\ref{lemma14}. Moreover, from \eqref{e91} we deduce \eqref{e92}.
Inequalities \eqref{e93}-\eqref{e94} can be shown similarly.
\hfill $\square$

\smallskip
\medskip

Note that in view of \eqref{e95}, for all $n\in \ene$ we can
define the unit tangent field
\[ \nu^n:=\frac{\nabla u^n}{|\nabla u^n|}\quad \mathcal V-\textrm{a.e. in}\;\;M\,,\]
and the dual unit cotangent field $\tilde \nu^n\,\,\big(\,i.e.\,\;
\tilde \nu^n(\nu^n)=1\,$ a.e. in $M\,\big)$.

\medskip

We have the following auxiliary identity.
\begin{lemma}\label{lemma17}
Let $\phi\in C^2_c(M; \re^+), \,U\subset\subset\{\phi>0\},\,
Y\equiv (Y^1,\ldots, Y^N)\in C^1_c(U; TM)$. There holds:
\begin{equation}\label{e140}
\begin{split}
-\e \Delta u^n \langle\nabla u^n, Y \rangle\,=\, \frac
{\e}2\di\left(Y|\nabla u^n|^2 \right)-\frac{\e}2|\nabla u^n|^2\cd Y : I\\
-\e \di(\nabla u^n\langle\nabla u^n, Y \rangle)+\e \cd Y:\nabla
u^n\otimes du^n\,.\hspace{.6 cm}
\end{split}
\end{equation}
\end{lemma}

\noindent{\it Proof \,.} Write $u$ instead of $u^n$ for brevity. Take any $p\in M$ and fix an orthonormal
frame $\{E_i\}_{i=1,\ldots, N}$ around $p$. We have $\nabla u=
\sum_{i=1}^N E_i u E_i\quad \textrm{around}\;p \,;$ furthermore,
if $Z=\sum_{i=1}^N Z^i E_i,$ then at $p$ there holds $\quad \di
Z=\sum_{i=1}^N E_i Z^i\,$ and
\[\Delta u E_iu =\Big(\sum_{j=1}^N E_j E_j u\Big) E_i u=\sum_{j=1}^N\Big[E_j(E_j uE_i u)-E_j uE_j(E_i u)\Big]\]
\[=\sum_{j=1}^N E_j(E_juE_i u)-\frac 12 E_i\sum_{j=1}^N(E_j u)^2;\]
here equality $[E_i, E_j](p)=0$ has been used. Thus,
\begin{equation}\label{e141}
\begin{split}
-\e \Delta u\langle \nabla u, Y\rangle=-\e \Delta u\sum_{i=1}^N
E_i u Y^i \hspace{1.2 cm}\\
=-\e \sum_{i=1}^N Y^i \sum_{j=1}^N E_j(E_j u E_i
u)+\frac{\e}2\sum_{i=1}^N Y^iE_i|\nabla u|^2\hspace{.3 cm}\\
=-\e \sum_{j,i=1}^N Y^i E_j(E_ju E_i u)+\frac{\e}2
\sum_{i=1}^N\big[E_i(Y^i|\nabla u|^2)-E_iY^i|\nabla u|^2\big]\\
=\e \sum_{j,i=1}^N\big[ - E_j(Y^i E_j uE_i u)+ E_j Y^iE_j uE_i
u\big]\hspace{.7 cm}\\
+\frac{\e}2 \di(Y|\nabla u|^2)-\frac{\e}2
|\nabla u|^2\di Y\hspace{1.3 cm}\\
= -\e \sum_{j=1}^N E_j(E_j u\langle\nabla u, Y \rangle)+\e
\sum_{i=1}^N \langle \nabla Y^i,\nabla u  \rangle E_i u\hspace{.5 cm}\\
+\frac{\e}2 \di(Y|\nabla u|^2)-\frac{\e}2 |\nabla u|^2\di
Y\,.\hspace{1.3 cm}
\end{split}
\end{equation}
It is easily seen that
\begin{equation}\label{e142}
\cd Y: I = \di Y\,,
\end{equation}
\begin{equation}\label{e143}
\quad \cd Y:\nabla u\otimes du = \sum_{j,i=1}^N E_i Y^j E_j uE_i
u=\sum_{i=1}^N \langle \nabla Y^i, \nabla u\rangle E_i u\,.
\end{equation}
From \eqref{e141}-\eqref{e143} we get \eqref{e140}. \hfill
$\square$

\medskip
\smallskip

The following representation formula for $\del V^n$ holds.

\begin{lemma}\label{18} Let assumptions of Lemma
\ref{lemma17} be satisfied. There holds:
\begin{equation}\label{e96}
\begin{split}
\del V^n(Y)=\int_{M}\nu^n\otimes\tilde \nu^n:\cd Y d\xi^n \hspace{1 cm}\\
-\int_{M} \e_n \langle Y, \nabla u^n\rangle\left(-\Delta u^n
+\frac{f(u^n)}{\e_n^2}\right)d \mathcal V(x)\,.
\end{split}
\end{equation}
\end{lemma}

\noindent{\it Proof\,.}
Defining the {\it stress} tensor by
\[T:= \left\{\frac{\e_n}2|\nabla u^n|^2+\frac 1{\e_n}F(u^n)\right\} I-\e_n \nabla u^n \otimes d u^n\,,\]
we have
\begin{equation}\label{e173}
T= \frac{\e_n}{2}|\nabla u^n|^2\left(I- 2\nu^n \otimes\tilde
\nu^n\right) + \frac{F(u^n)}{\e_n} I\,.
\end{equation}

Integrating by parts, by \eqref{e140} and
\eqref{e173},
\begin{equation}\label{e111}
\begin{split}
\int_{M}\e_n \left(-\Delta u^n +
\frac{f(u^n)}{\e^2_n}\right)\langle \nabla u^n, Y \rangle d
\mathcal V=-\int_{M} T : \cd Y d\mathcal V \\
=-\int_{M}\left(\frac{\e_n}{2}|\nabla
u^n|^2+\frac{F(u^n)}{\e_n}\right)(I- \nu^n\otimes\tilde \nu^n):
\cd Y
d\mathcal V\hspace{.5 cm}\\
+\int_{M}\left(-\frac{F(u^n)}{\e_n}+\frac{\e_n}2|\nabla
u^n|^2\right)\nu^n\otimes\tilde \nu^n:\cd Y d\mathcal V \hspace{.3 cm}\\
=-\int_{M}(I-\nu^n\otimes\tilde \nu^n): \cd Y d\mu^n+\int_{M}
\nu^n\otimes\tilde \nu^n: \cd Y d\xi^n\,. \hspace{.2 cm}
\end{split}
\end{equation}
Since
\[\del V^n(Y)=\int_{M} \cd Y : S \, dV(x,S)\hspace{1.6 cm}\\
= \int_{M} \cd Y : (I-\nu^n\otimes\nu^n) d\mu^{n}\,,\] from
\eqref{e111} we obtain \eqref{e96}. \hfill $\square$

\medskip

Now we are ready to prove Proposition \ref{prop7}.

\smallskip

\noindent{\it Proof of Proposition \ref{prop7}\,.\,\,} Keep the
same notation as above. By compactness theorem for Radon measures
with locally equibounded masses, there exist a subsequence of
$\{V^{n_k}\}\subset\{V^n\}$ and $\tilde V\in \mathbb V_{N-1}(M)$
such that $V^n\to \tilde V$ as $n\to \infty,$ as varifolds. Let us
write $\{V^{n_k}\}\equiv\{V^n\}$.

\medskip

\noindent {\bf Claim 1.} We have that $\tilde
V\lfloor\{\phi>0\}\in \mathcal R \mathbb V_{N-1}(\{\phi>0\})$;
moreover, $\|\tilde V\|=\mu$ and $\mu$ is a $(N-1)-$rectifiable
Radon measure on $M$.

\smallskip

In fact, for $U\subset\subset \{\phi>0\}$ and $Y\in C^1_c(U; TM)$,
sending $n\to \infty$ in \eqref{e96}, in view of hypothesis
$(ii)$, we get
\begin{equation}\label{e97}
\begin{split}
\del \tilde V(Y)=\int_{M} \cd Y : S\, d\tilde V(x,S)\\
=-\lim_{n\to \infty}\int_{M}\e_n \langle Y, \nabla
u^n\rangle\left(-\Delta u^n +\frac{f(u^n)}{\e^2_n}\right) d
\mathcal V.
\end{split}
\end{equation}
Furthermore, by hypothesis $(iii)$ and \eqref{e94} we obtain
\[|\del \tilde V(Y)|\le \|Y\|_\infty\limsup_{n\to \infty}\int_U \e_n |\nabla u^n|\left|-\Delta u^n+\frac{f(u^n)}{\e_n^2}\right| d \mathcal V\]
\[\le \|Y\|_\infty\limsup_{n\to \infty}\int_U\left\{\frac{\e_n}2\frac{|\nabla u^n|^2}{\phi} +\e_n\phi\left(-\Delta u^n+\frac{f(u^n)}{\e_n^2}
\right)^2\right\}d \mathcal V \]
\[\le \|Y\|_\infty \limsup_{n\to \infty}\Big[2 C(\phi, U)\mu^n\big(\{\phi>0\}\big)+ 2 \check C + 4
C_1(\phi)\mu^n\big(\{\phi>0\}\big)\Big]\]
\[=\|Y\|_\infty\Big[C(\phi, U)\mu\big(\{\phi>0\}\big)+ 2 \check C\Big]\,.\]
Hence
\begin{equation}\label{e98}
|\d \tilde V(Y)|\le C(\phi, U, \mu,\check C)\|Y\|_\infty.
\end{equation}
By \eqref{e98}, $|\del \tilde V|\lfloor\{\phi>0\}$ is a Radon
measure on $\{\phi>0\}$. By Corollary $2)$ of the rectifiability
theorem in \cite[pag. 450]{All} and hypothesis $(iv)$, Claim 1
follows. So, $(a)$ has been verified.

\smallskip

In view of rectifiability, the varifold $V\equiv \tilde V$ is
uniquely determined by $\mu$, independently of the subsequence;
thus, the all sequence $\big\{V^n\lfloor\{\phi>0\} \big\}$
converges to $V$ as varifolds. So, $(b)$ follows. From \eqref{e97}
we get $(c)$.

\smallskip

It remains to prove $(d)$. To this aim take $\psi\in
C^2_c\big(\{\phi>0\}; \re^+\big)$ with $\sqrt \psi\in
C^1(\{\phi>0\}; \re^+)$. Since $\mu$ is rectifiable, in view of
Lemma \ref{lemma15} we have:
\[\Big(\int_{\{\phi>0\}} \psi H^2 d\mu \Big)^{1/2} =
\sup\left\{ \langle \sqrt \psi Y,  \overrightarrow{H} \rangle :
Y\in C^\infty_c(\{\phi>0\}; TM), \|Y\|_{L^2(\mu)}\le 1 \right\}\,.
\]

By \eqref{e96}, with $Y$ replaced by $\sqrt \psi Y$, and
assumption $(ii)$, we get:
\[\int_{\suppo \phi} \sqrt \psi \langle Y, \overrightarrow{H}\rangle d\mu = - \d V(\sqrt \psi Y)=-\lim_{n\to\infty} \d V^n(\sqrt \psi Y)\]
\[= \lim_{n\to\infty} \int_{\suppo \phi}\e_n \sqrt \psi \langle Y, \nabla u^n \rangle\left(-\Delta u^n +\frac{f(u^n)}{\e_n^2}\right)d \mathcal V\]
\[-\lim_{n\to \infty} \int_{\suppo \phi} \nu^n \otimes \tilde{\nu}^n : \cd (\sqrt \psi Y) d\xi^n\]
\[\le \limsup_{n\to\infty} \left(\int_{\suppo \phi}\e_n |\nabla u^n|^2 |Y|^2 d \mathcal V\right)^{1/2}\left[\int_{\suppo \phi}\e_n\psi\left(-\Delta u^n+
\frac{f(u^n)}{\e_n^2} \right)^2 \right]^{1/2}d\mathcal V \]
\[\le \limsup_{n\to \infty} \left(\int_{\suppo \phi}|Y|^2d\mu^n \right)^{1/2}\limsup_{n\to \infty}\left[\int_{\suppo \phi}\e_n\psi
\left(-\Delta u^n +\frac{f(u^n)}{\e_n^2}\right)^2 d \mathcal V
\right]^{1/2}\]
\[=\|Y\|_{L^2(\mu)}\limsup_{n\to \infty}\left[\int_{\suppo \phi}\e_n \psi\left(-\Delta u^n +\frac{f(u^n)}{\e_n^2}\right)^2 d \mathcal V\right]^{1/2}\,.\]
Fixing $\psi=\psi_k$ and taking the $\sup$ on $Y$, this implies
\begin{equation}\label{e99}
\int_{\suppo \phi}\psi_k H^2 d\mu\le \limsup_{n\to \infty}
\int_{\suppo \phi} \e_n\psi_k\left(-\Delta u^n
+\frac{f(u^n)}{\e_n^2}\right)^2d \mathcal V \,
\end{equation}
for all \,$k\in \ene$, where $\{\psi_k\}_{k\in \ene}\subset
C^2_c(\{\phi>0\})$ with $\sqrt{\psi_k}\in
C^1(\{\phi>0\}),\,\psi_k\le \psi_{k+1}$ for every $k\in \ene,
\psi_k\to \phi$ in $L^1(\{\phi>0\})$. Letting $k\to \infty$ in
\eqref{e99}, in view of the monotone convergence theorem we have:
\begin{equation}\label{e100}
\begin{split}
\int_{\suppo \phi} \phi H^2 d\mu \le \limsup_{n\to \infty}
\int_{\suppo \phi} \e_n \phi \left(-\Delta u^n
+\frac{f(u^n)}{\e_n^2} \right)^2d \mathcal V \\\le C(\phi, \mu,
\check C)\,;\hspace{3 cm}
\end{split}
\end{equation}
in the last inequality \eqref{e94} and hypothesis $(iii)$   have
been used.

\medskip
\smallskip

\noindent {\bf Claim 2.} The following equality holds:
\begin{equation}\label{e101}
\begin{split}
\lim_{n\to \infty} \int_{\suppo \psi} \e_n \langle \nabla \psi,
\nabla u^n \rangle\left(-\Delta u^n +\frac{f(u^n)}{\e_n^2}\right)d
\mathcal V\\= \int_{\suppo \psi} \langle\nabla \psi, S^{\perp}
(\overrightarrow{H})\rangle d\mu\,, \hspace{1.5 cm}
\end{split}
\end{equation}
for any $\psi\in C^2_c\big(\{\phi>0\};TM\big)$\,; here
$S=S(x)=T_x\, \mu$\,.

\smallskip

To prove Claim 2, note that since $\mu$ is rectifiable, by Lemma
\ref{lemma15} for each $\d>0$ we can select $Y\in
C^1_c\big(\{\phi>0\}\,; TM \big)$ such that
\begin{equation}\label{e101a}
\int_{\suppo \psi} |Y(x)- S^{\perp}(\nabla\psi(x))|^2 d\mu\le
\d^2\,.
\end{equation}
From \eqref{e96} we obtain:
\begin{equation}\label{e102}
\int_{\suppo \psi} \langle\nabla \psi, S^{\perp}
(\overrightarrow{H})\rangle\, d\mu =\int_{\suppo \psi} \langle
S^{\perp}(\nabla \psi),  \overrightarrow{H}\rangle\,
d\mu\,=\,\sum_{i=1}^6\mathcal A_i\,,
\end{equation}
where
\[\mathcal A_1:= \int_{\suppo \psi} \Big(S^{\perp}(\nabla \psi)  - Y\Big)
\overrightarrow{H}\, d\mu\,,\]
\[\mathcal A_2:=-\d V(Y) +\d V^n(Y)\]
\[\mathcal A_3:=-\d V^n(Y) -\int_{\suppo \psi} \e_n \langle\nabla u^n, Y
\rangle\left(-\Delta u^n+\frac{f(u^n)}{\e_n^2}\right)d \mathcal
V\]
\[\mathcal A_4:=\int_{\suppo \psi} \e_n \big\langle Y-\langle \nabla \psi,
\nu^n\rangle \nu^n, \nabla u^n \big\rangle \left( -\Delta
u^n+\frac{f(u^n)}{\e_n^2}\right)d \mathcal V\]
\[\mathcal A_5:=\int_{\suppo \psi}\e_n
\big\langle \langle \nabla \psi, \nu^n\rangle \nu^n-\nabla \psi,
\nabla u^n \big\rangle \left( -\Delta
u^n+\frac{f(u^n)}{\e_n^2}\right)d \mathcal V\]
\[\mathcal
A_6:=\int_{\suppo \psi} \e_n \langle \nabla \psi, \nabla
u^n\rangle \left(-\Delta u^n+\frac{f(u^n)}{\e_n^2}\right)d
\mathcal V\,.\] Let us estimate $|\mathcal A_i|$ for $i=1,\ldots,
6\,.$

From \eqref{e100} and \eqref{e101a} we get
\begin{equation}\label{e103}
|\mathcal A_1|\le \d \left(\int_{\suppo \psi} H^2 d\mu
\right)^{1/2}\le C(\psi, \phi)\, \d\,.
\end{equation}
In view of $(b)$, we get
\begin{equation}\label{e104}
\lim_{n\to \infty} \mathcal A_2\,=\,0\,.
\end{equation}
By \eqref{e97},
\begin{equation}\label{e105}
\lim_{n\to \infty} \mathcal A_3\,=\lim_{n\to \infty}\int_{\suppo
\psi} \nu^n\otimes \tilde{\nu}^n : \cd Y\, d\xi^n\,=\,0\,.
\end{equation}
Moreover, by \eqref{e94}, hypothesis $(iii)$ and H$\ddot{o}$lder
inequality,
\[
|\mathcal A_4|\le \left(\e_n\int_{\suppo \psi}|\nabla
u^n|^2\,\Big| Y- \langle \nabla \psi, \nu^n\rangle\nu^n\Big|^2 d
\mathcal V\right)^{1/2}\]
\[\cdot \left[\sup_{\suppo \psi}\frac 1{\phi}\int_{\suppo
\psi}\e_n\phi\left(-\Delta u^n +\frac{f(u^n)}{\e_n^2}\right)^2 d
\mathcal V \right]^{\frac 1 2} \]
\[ \le \left(2\int_{\suppo \psi}\Big|Y(x)-S^{\perp}\big(\nabla
\psi(x)\big)\Big|^2 dV^n(x,S) \right)^{1/2} [C(\suppo \psi,
\phi)]^{1/2}
\]
\[ \cdot \left[-2 \mathfrak B^{n}(u^n, \phi) + 4
C_1(\phi)\int_{\suppo \psi}\frac{\e_n}2|\nabla u^n|^2d \mathcal V
\right]^{1/2}\]
\[ \le \left(2\int_{M}\Big|Y- S^{\perp}(\nabla \psi)\Big|^2
dV^n \right)^{1/2} [C(\suppo \psi, \phi)]^{1/2}[C(\phi,
\mu^{\e_n}_t(\{\phi>0\}), \check C)]^{1/2}\,.
\]
Thus,
\begin{equation}\label{e106}
\begin{split}
\limsup_{n\to \infty} |\mathcal A_4|\le C(\phi,\psi)\big(
2\int_{\suppo \psi}|Y-S^{\perp}(\nabla \psi) |^2 dV^n
\big)^{1/2}\\ \le C(\phi, \psi, \mu, \check C)\d\,. \hspace{2.3
cm}
\end{split}
\end{equation}
Note that by definition of $\nu^n$, $\mathcal I_5\equiv 0$, hence
from \eqref{e102}-\eqref{e106} we deduce
\[\left| \int_{\suppo \psi} \e_n \langle \nabla \psi, \nabla
u^n\rangle \left(-\Delta u^n+\frac{f(u^n)}{\e_n^2}\right)d
\mathcal V - \int_{\suppo \psi} \langle S^{\perp}(\nabla \psi),
\overrightarrow{H}\rangle\, d\mu_t\right|\]\[\le C(\psi,
\phi)\d\,.\] Letting $\d\to 0^+,$ we obtain \eqref{e101}.

\medskip
\smallskip

\noindent{\bf Claim 3:} The limit \eqref{e101} remains true with
$\psi$ replaced by $\phi$.

\smallskip

Too see this, for any $\d>0$ we take $\psi\in
C^2_c\big(\{\phi>0\};\re^+\big)$ such that $\psi<\phi,\,
\|\phi-\psi\|_{C^2}<\d\,.$ We write:
\begin{equation}\label{e107}
\begin{split}
\limsup_{n\to \infty}\Big| \int_{\suppo \phi} \e_n \langle \nabla
\phi, \nabla u^n \rangle\left(-\Delta u^n
+\frac{f(u^n)}{\e_n^2}\right)d \mathcal V\\ - \int_{\suppo \phi}
\langle S^{\perp}(\nabla \phi), \overrightarrow{H}\rangle
d\mu\Big| \le \limsup_{n\to\infty} |\mathcal I_1+\mathcal
I_2+\mathcal I_3|,\hspace{.2 cm}
\end{split}
\end{equation}
where
\[\mathcal I_1:= - \int_{\suppo \psi}\langle  S^{\perp}(\nabla \psi),\overrightarrow{H} \rangle d\mu
+ \lim_{n\to \infty} \int_{\suppo \psi} \e_n \langle \nabla \psi,
\nabla u^n \rangle\left(-\Delta u^n +\frac{f(u^n)}{\e_n^2}\right)d
\mathcal V\,,\]
\[\mathcal I_2:=\lim_{n\to \infty}\int_{\suppo \phi}\e_n \langle\nabla \phi-\nabla \psi, \nabla u^n\rangle
\left(-\Delta u^n +\frac{f(u^n)}{\e_n^2}\right)d \mathcal V\,,\]
\[\mathcal I_3:= \int_{\suppo \phi} \langle S^{\perp}(\nabla \phi-\nabla \psi),\overrightarrow{H}\rangle d\mu\,.\]
By \eqref{e101},
\begin{equation}\label{e108}
\limsup_{n\to\infty}|\mathcal I_1| = 0\,.
\end{equation}
By \eqref{e100}, H\"older inequality and Lemma
\ref{lemma14},
\begin{equation}\label{e109}
\begin{split}
\limsup_{n\to\infty}|\mathcal I_2|\le \limsup_{n\to
\infty}\left(\int_{\{\phi>0\}}\frac{|\nabla\phi-\nabla
\psi|^2}{\phi-\psi}|\nabla
u^n|^2d \mathcal V \right)^{1/2}\hspace{.2 cm}\\
\cdot \left[\int_{\{\phi>0\}} \e_n (\phi-\psi)\left(-\Delta u^n +\frac{f(u^n)}{\e_n^2}\right)^2 d \mathcal V\right]^{1/2} \hspace{.8 cm}\\
\le
[2\mu(\{\phi>0\})\sup_{\{\phi>0\}}|\Hess(\phi-\psi)|]^{1/2}C(\phi,\mu,\check
C)^{1/2}\le C(\phi, \mu,\check C)\d\,.
\end{split}
\end{equation}
Furthermore, it is easily checked that from \eqref{e100} one also
has
\begin{equation}\label{e110}
|\mathcal I_3|\le C(\phi)\d\,.
\end{equation}
From \eqref{e107}-\eqref{e110}, letting $\d\to 0^+$, and combining
\eqref{e99} and \eqref{e100} with \eqref{e171} and \eqref{e172},
we get Claim 3 and $(d)$. This completes the proof. \hfill
$\square$

\subsection{Proof of Theorem \ref{thmdb}}

Finally we can prove the main result of the paper.

\medskip
\noindent{\it Proof of Theorem \ref{thmdb}\,.\,\,}

For any $\phi\in C^2_c(M, \re^+)$, from the semidecreasing property in Lemma
\ref{lemma7} we have that
\begin{equation}\label{e112b}
\ud_{ t_0}\mu_t(\phi)>-\infty
\end{equation}
for a.e. $t_0>0$, the upper derivative being clearly an ordinary derivative. Now, fix any $t_0 >0$ such that
\eqref{e112b} is satisfied, otherwise there is nothing to be proven. Set
\begin{equation}\label{e112a}
-\infty<\tilde D\equiv \ud_{t_0}\mu_t(\phi)\,.
\end{equation}
Thus, there exists sequences $\{\d_k\}\subset (0,1),\, \d_k\to 0$
and $\{t_k\}\subset (0,\infty), t_k\to t_0$ as $k\to \infty$ such
that
\[\tilde D - \d_k \le \frac{\mu_{t_k}(\phi)-\mu_{t_0}(\phi)}{t_k-t_0};\]
we may assume that $t_k>t_0$ for all $k\in \ene.$

Since $\mu_t^{\e_n}\to \mu_t$, we can find a sequence
$\{r_k\}\subset (0,\infty), r_k\to \infty$ such that for all
$k\in\ene$
\begin{equation}
\tilde D - 2\d_k \le
\frac{\mu_{t_k}^{\e_{r_k}}(\phi)-\mu_{t_0}^{\e_{r_k}}(\phi)}{t_k-t_0}=\frac
1 {t_k-t_0}\int_{t_0}^{t_k}\frac{d}{dt}
\mu_t^{\e_{r_k}}(\phi)dt\,.
\end{equation}
Note that, for $d|\xi^{\e_n}|\to 0$ in $M\times [0,\infty)$, we
can increase $r_k$ so that
\begin{equation}\label{e112}
\int_{t_0}^{t_k}\int_{\{\phi>0\}} d |\xi^{\e_{r_k}}|\, \le \d_k^2
(t_k - t_0)\,.
\end{equation}
By the proof of Lemma \ref{lemma7}, there exists $C_1=C_1(\phi)>0$ such that
\[\frac d{dt}\mu^{\e_n}_t(\phi)\le C_1(\phi)\quad \textrm{for all}\;\; n\in \ene, 0<t<t_0+1\,.\]
For any $k\in \ene$, define
\[\mathcal Z_k:=\left\{t\in [t_0, t_k]\,:\,\frac{d}{dt}\mu^{\e_{r_k}}(\phi)\ge \tilde D-3\d_k\right\}\,.\]
We have
\[\tilde D-2\d_k \le \frac 1{t_k -t_0}\int_{[t_0, t_k]\setminus \mathcal Z_k}(\tilde D-3\d_k)dt+\frac 1{t_k-t_0}\int_{\mathcal Z_k} C_1(\phi) dt.\]
So,
\[meas(\mathcal Z_k)\ge \frac{\d_k(t_k-t_0)}{C_1(\phi)-\tilde D+ 3\d_k}\ge \frac{\d_k(t_k-t_0)}{2(C_1(\phi)-\tilde D)},\]
for $k\in \ene$ big enough. By \eqref{e112},
\[meas(\mathcal Z_k)\inf_{t\in \mathcal Z_k} \big|\xi_{t}^{\e_{r_k}}\big|\big(\{\phi>0\}\le \d_k^2 (t_k-t_0)\big)\,.\]
Hence, due to \eqref{e175} and \eqref{e173}, we can construct a
sequence $\{s_k\}\subset \mathcal Z_k$ such that
\begin{equation}\label{e113}
\tilde D-3\d_k\le \frac{d}{dt}\mu_t^{\e_{r_k}}(\phi)|_{t=s_k} =
\mathcal B^{\e_{r_k}}(u^{\e_{r_k}}(\cdot, s_k), \phi)\,,
\end{equation}
and
\begin{equation}\label{e114}
\left|\xi_{s_k}^{\e_{r_k}}\right|\big(\{\phi>0\}\big)\le 2 (C_1-
\tilde D)\d_n\,.
\end{equation}

By hypothesis \eqref{e170} with $K=\suppo \phi$ and standard
compactness results, there exists a subsequence of
$\{\mu^{\e_{r_k}}_{s_k}\}$, which converges to $\tilde\mu$, for
some Radon measure on $M$ $\tilde \mu$. By Lemma \ref{lemma7} it
is possible to show that (see \cite{Ilm2}, Section 7.6)
\begin{equation}\label{e115}
\tilde\mu\lfloor\{\phi>0\}=\mu_{t_0}\lfloor\{\phi>0\},
\end{equation}
hence $\mathfrak B(\tilde\mu,\phi)=\mathfrak B(\mu_{t_0},\phi)$\,.

\smallskip
Corollary \ref{corhauss} combined with \eqref{e113}-\eqref{e115}
implies that hypotheses $(i)-(iv)$ of Proposition \ref{prop7} are
satisfied with $\{u^n\}$ replaced by $\{u^k\}$, where $u^k\equiv
u^{\e_{r_k}}(\cdot, s_k)$. By \eqref{e112a},
\eqref{e113}-\eqref{e115}, due to Proposition \ref{prop7}-$(a)$,
we see that $\mu_{t_0}$ is locally $(N-1)-$rectifiable (varying $\phi=\phi_i \in C^2_c(M,\mathbb{R}^+)$ in a coutable set of functions such that $\cup \{\phi_i>0\}=M$). Moreover, in view
of Proposition \ref{prop7}-$(d)$, we obtain
\[\ud_{t_0} \mu_t(\phi)\, \le\, \mathfrak B(\mu_{t_0},\phi).\]
This completes the proof. \hfill $\square$

\end{document}